\documentclass{amsart}
\usepackage{times,mathrsfs,array,amssymb,pb-diagram,epsfig}

\newcounter{lemma}
\newtheorem{Theorem}{Theorem}
\newtheorem{Lemma}[lemma]{Lemma}
\newtheorem{Corollary}[lemma]{Corollary}
\newtheorem{Proposition}[lemma]{Proposition}
\newtheorem{theorem}{Theorem}

\theoremstyle{definition}
\newtheorem{Example}[lemma]{Example}

\newtheorem{Remark}[lemma]{Remark}

\def\C{\mathbb C}
\def\H{\mathbb H}
\def\O{\mathcal O}
\def\P{\mathbb P}
\def\Q{\mathbb Q}
\def\R{\mathbb R}
\def\Z{\mathbb Z}
\def\fp{\mathfrak p}

\def\Im{\mathrm{Im}\,}
\def\FL#1{\left\lfloor #1\right\rfloor}

\def\div{\operatorname{div}}
\def\mod{\  \mathrm{mod}\ }

\def\sgn{\mathrm{sgn\,}}
\def\Im{\mathrm{Im\,}}
\def\n{\mathit{n}}
\def\tr{\mathit{tr}}

\def\gen#1{\langle #1\rangle}
\def\JS#1#2{\left(\frac{#1}{#2}\right)}
\def\GL{\mathrm{GL}}
\def\Gr{\mathrm{Gr}}
\def\SL{\mathrm{SL}}
\def\CM{\mathrm{CM}}
\def\M#1#2#3#4{\begin{pmatrix}#1&#2\\#3&#4\end{pmatrix}}
\def\SM#1#2#3#4{\left(\begin{smallmatrix}#1&#2\\#3&#4\end{smallmatrix}
  \right)}
\def\wt{\widetilde}

\begin{document}

\title{Special values of hypergeometric functions and periods of CM
  elliptic curves}

\author{Yifan Yang}
\address{Department of Applied Mathematics, National Chiao Tung
  University and National Center for Theoretical Sciences, Hsinchu,
  Taiwan 300}  
\email{yfyang@math.nctu.edu.tw}
\date{\today}
\subjclass[2000]{primary 11F12; secondary 11G15, 11G18, 33C05}
\thanks{The author was partially supported by Grant
  102-2115-M-009－001-MY4 of the National Science Council, Taiwan (R.O.C.).}

\begin{abstract} Let $X_0^6(1)/W_6$ be the Atkin-Lehner quotient of
  the Shimura curve $X_0^6(1)$ associated to a maximal order in an
  indefinite quaternion algebra of discriminant $6$ over $\Q$. By
  realizing modular forms on $X_0^6(1)/W_6$ in two ways, one in terms
  of hypergeometric functions and the other in terms of Borcherds
  forms, and using Schofer's formula for values of Borcherds forms at
  CM-points, we obtain special values of certain hypergeometric
  functions in terms of periods of elliptic curves over $\overline\Q$
  with complex multiplication.
\end{abstract}

\thanks{The author would like to thank Srinath Baba, H\aa kan Granath,
  and John Voight for many fruitful discussions.}

\maketitle

\begin{section}{Introduction}

Let $X_0^D(N)$ be the Shimura curve associated to an Eichler order of
level $N$ in an indefinite quaternion algebra of discriminant $D$ over
$\Q$. When $D=1$, the Shimura curve $X_0^1(N)$ is just the classical
modular curve $X_0(N)$ and there are many different constructions of
modular forms on $X_0(N)$ in literature, such as Eisenstein series,
Dedekind eta functions, Poincare series, theta series, and etc. These
explicit constructions provide practical tools for solving problems
related to classical modular curves. On the other hand, when $D\neq
1$, because of the lack of cusps, most of the methods for classical
modular curves cannot possibly be extended to the case of general
Shimura curves. As a result, even some of the most fundamental
problems about Shimura curves, such as finding equations of Shimura
curves, computing Hecke operators on explicitly given modular forms,
and etc., are not easy to answer. However, in recent years, there have
been two realizations of modular forms on Shimura curve emerging in
literature and some progress toward the study of Shimura curves has
already been made using these two methods.

The first method was due to the author of the present paper. In
\cite{Yang-Schwarzian}, we first observed that when a Shimura curve $X$
has genus $0$, all modular forms on $X$ can be expressed in terms
of solutions of the Schwarzian differential equation associated to a
Hauptmodul of $X$. Then by utilizing the Jacquet-Langlands
correspondence and explicit covers between Shimura curves, we
devised a method to compute Hecke operators with respect to the
explicitly given basis of modular forms. As applications of this
computation of Hecke operators, we computed modular equations for
Shimura curves, which can be regarded as equations for Shimura curves
associated to Eichler orders of higher levels, in \cite{Yang-ModEqs}
and obtained Ramanujan-type identities for Shimura curves in
\cite{Yang-Ramanujan}. In addition, since some Schwarzian differential
equations are essentially hypergeometric differential equations, this
realization of modular forms yields many beautiful identities among
hypergeometric functions. This is discussed in \cite{Tu,Tu-Yang}.

The second method is to realize meromorphic modular forms with
divisors supported on CM-points as Borcherds forms associated to the
lattice formed by the elements of trace zero in an Eichler order.
Borcherds forms themselves are not easy to work with. What makes
Borcherds forms useful in practice is Schofer's formula \cite{Schofer}
for norms of (generalized) singular moduli of Borcherds forms, that
is, norms of values of Borcherds forms at CM-points. Schofer's formula
is based on an earlier work of Kudla \cite{Kudla}, and the evaluation of
derivatives of Fourier coefficients of Eisenstein series uses
works of Kudla, Rapoport, and Yang
\cite{Kudla-Rapoport-Yang,Kudla-Yang,Yang-density}. An immediate
consequence of Schofer's formula is a necessary condition for primes
that can appear in the prime factorization of the norm of the
difference of two singular moduli of different discriminants, which is
analogous to Gross and Zagier's work \cite{Gross-Zagier} for the case
of the classical modular curve $X_0(1)$. Also, Errthum \cite{Errthum}
used Schofer's formula to determine singular moduli of $X_0^6(1)/W_6$
and $X_0^{10}(1)/W_{10}$, where $W_D$ denotes the group of all
Atkin-Lehner involutions on $X_0^D(1)$, verifying Elkies' numerical
computation \cite{Elkies}. (However, we remark that Schofer's formula
needs a slight correction when the Borcherds forms have nonzero weights. 
See Section \ref{section: Schofer} below.)

The realization of modular forms on Shimura curves in
\cite{Yang-Schwarzian} is completely analytic, while Schofer's formula
for singular moduli of Borcherds forms is more arithmetic in nature.
(For example, the primary motivation of
\cite{Kudla,Kudla-Rapoport-Yang, Kudla-Rapoport-Yang2,Kudla-Yang} was
to obtain arithmetic Siegel-Weil formulas realizing generating series
from arithmetic geometry as modular forms.) It is an interesting
problem to see what results we can obtain by combining the two
approaches. This is the main motivation of the present work.

In this paper, we will consider the Shimura curve $X=X_0^6(1)/W_6$.
From \cite{Yang-Schwarzian}, we know that every holomorphic modular
form on $X$ can be expressed in terms of hypergeometric functions.
Now according to \cite[Theorem 7.1]{Shimura-periods} and \cite[Theorem
1.2 and (1.4) of Chapter 3]{Yoshida}, if $t(\tau)$ is a modular
function on $X$ that takes algebraic values at all CM-points, then the
value of $t'(\tau)$ at a CM-point of discriminant $d$
is an algebraic multiple of the square of
$$
  \omega_d=e^{L'(0,\chi_{d_0})/2L(0,\chi_{d_0})}=\frac1{\sqrt{|d_0|}}
  \prod_{a=1}^{|d_0|-1}\Gamma\left(\frac a{|d_0|}\right)^{\chi_{d_0}(a)\mu_{d_0}/4h_{d_0}},
$$
where $d_0$ is the discriminant of the field $\Q(\sqrt d)$,
$\chi_{d_0}$ is the Kronecker character associated to 
$\Q(\sqrt d)$, $\mu_{d_0}$ is the number of roots of unity in $\Q(\sqrt
d)$, and $h_{d_0}$ is the class number of $\Q(\sqrt d)$. (See
\cite[Theorem 7]{Baba-Granath} for some examples.) The significance of
these numbers $\omega_d$ is that periods of any elliptic curve over
$\overline\Q$ with CM by $\Q(\sqrt d)$ lie in
$\sqrt\pi\omega_d\cdot\overline\Q$. (See \cite{Gross,Chowla-Selberg}.)
In other words, the values of certain hypergeometric functions at
singular moduli can be expressed in terms of periods of CM elliptic
curves over $\overline\Q$.

\begin{Theorem} \label{theorem: 0}
  Let $s(\tau)$ be the Hauptmodul of $X_0^6(1)/W_6$ that takes values
  $0$, $1$, and $\infty$ at the CM-points of discriminants $-4$,
  $-24$, and $-3$, respectively. Let $\tau_d$ be a CM-point of
  discriminant $d$ such that $|s(\tau_d)|<1$. Then
  $$
    {}_2F_1\left(\frac1{24},\frac5{24};\frac34;s(\tau_d)\right)
    \in\frac{\omega_{d}}{\omega_{-4}}\cdot\overline\Q, \quad
    {}_3F_2\left(\frac13,\frac12,\frac23;\frac34,\frac54;s(\tau_d)\right)
    \in \omega_{d}^2\cdot\overline\Q.
  $$
  Likewise, let $t(\tau)=1/s(\tau)$. If $\tau_d$ is a CM-point of
  discriminant $d$ such that $|t(\tau_d)|<1$, then
  $$
    {}_2F_1\left(\frac1{24},\frac7{24};\frac56;t(\tau_d)\right)
    \in\frac{\omega_{d}}{\omega_{-3}}\cdot\overline\Q, \quad
    {}_3F_2\left(\frac14,\frac12,\frac34;\frac56,\frac76;
    t(\tau_d)\right)\in\omega_{d}^2\cdot\overline\Q.
  $$
\end{Theorem}

The proof of the theorem will be given at the end of Section
\ref{section: Schwarzian}.

The parallel results in the cases of classical modular
curves can be described as follows. Let
$\lambda_1$ and $\lambda_2$ be a basis for a lattice $\Lambda$ in $\C$
with $\Im(\lambda_2/\lambda_1)>0$, and for positive even integers $k\ge 4$,
let
$$
  G_{k}(\Lambda)=\sum_{\lambda\in\Lambda,\lambda\neq 0}\frac1{\lambda^k}.
$$
Then Weierstrass's equation for the elliptic curve $\C/\Lambda$ over $\C$ is
$$
  y^2=4x^3-40G_4(\Lambda)x-140G_6(\Lambda).
$$
From the relations
$$
  G_4(\Lambda)=\frac1{45}\left(\frac\pi{\lambda_1}\right)^4E_4(\tau), \qquad
  G_6(\Lambda)=\frac2{945}\left(\frac\pi{\lambda_1}\right)^6
  E_6(\tau),
$$
where $\tau=\lambda_2/\lambda_1$ and $E_k$ are the normalized Eisentein
series of weight $k$, we immediately see that for
$\tau\in\Q(\sqrt d)\cap\H^+$, $\H^+=\{\tau\in\C:\Im\tau>0\}$,
$$
  E_k(\tau)\in \left(\frac{\Omega_d}\pi\right)^k\cdot\overline\Q,
$$
where $\Omega_d$ is any nonzero period of any elliptic curve over
$\overline\Q$ with CM by $\Q(\sqrt d)$. According to the
Chowla-Selberg formula \cite{Gross,Chowla-Selberg}, we may choose
$$
  \Omega_d=\sqrt\pi\prod_{a=1}^{|d|-1}\Gamma\left(\frac a{|d|}
  \right)^{\chi_d(a)\mu_d/4h_d}=\sqrt{\pi|d|}\omega_d.
$$
Now from the classical identity
$$
  E_4(\tau)={}_2F_1\left(\frac1{12},\frac5{12};1;\frac{1728}{j(\tau)}\right)^4,
$$
we conclude that if $\tau\in\Q(\sqrt d)\cap\H^+$, then
$$
  _2F_1\left(\frac1{12},\frac5{12};1;\frac{1728}{j(\tau)}\right)
  \in\frac{\Omega_d}\pi\cdot\overline\Q.
$$
For instance, for $\tau=i$, we have $j(i)=1728$, and Gauss' formula
for values of hypergeometric functions at $1$ and the 
multiplication formula for the Gamma function yield
\begin{equation*}
\begin{split}
  {}_2F_1\left(\frac1{12},\frac5{12};1;1\right)
  &=\frac{\Gamma(1/2)}{\Gamma(11/12)\Gamma(7/12)}
  =\frac{\sqrt\pi\Gamma(3/12)}{\Gamma(11/12)\Gamma(7/12)\Gamma(3/12)}
  \\
&=\frac{\sqrt\pi\Gamma(1/4)}{(2\pi)3^{1/2-3/4}\Gamma(3/4)}
 =\frac{3^{1/4}}2\frac{\Omega_{-4}}\pi.
\end{split}
\end{equation*}
For a fundamental discriminant $d<0$, one may use the
Chowla-Selberg formula
\cite[Page 110]{Chowla-Selberg}
$$
  \prod_{j=1}^{h_d}a_j^{-6}\Delta(\tau_j)=\frac{\omega_d^{12h_d}}{(2\pi)^{6h_d}},
$$
where the product runs through the complete set of reduced primitive
quadratic forms $a_jx^2+b_jxy+c_jy^2$ of discriminant $d$ with
$\tau_j=(-b_j+\sqrt d)/2a_j$, and its generalizations to
determine special values of hypergeometric functions.
See \cite{Archinard,Beukers-Wolfart,Chapman-Hart} for some examples.

Now to determine the precise values of the hypergeometric functions in
Theorem \ref{theorem: 0} at singular moduli, we shall realize the
modular forms involved as Borcherds forms. Then evaluating these
modular forms at CM-points using Schofer's formula, we obtain formulas
for special values of hypergeometric functions. The results in the
cases where there exists exactly one CM-point of fundamental
discriminant $d$ are given in the next theorem. In Section
\ref{section: additional}, we will work out an example to illustrate a
general technique to determine special values of the hypergeometric
functions when there are more than one CM-points of discriminant $d$.

\begin{Theorem} \label{theorem: main theorem}
The evaluations
\begin{equation*}
\begin{split}
  &{}_2F_1\left(\frac1{24},\frac5{24};\frac34;\frac MN\right)
  =A_1\frac{\omega_d}{\omega_{-4}}, \\
  &{}_3F_2\left(\frac13,\frac12,\frac23;\frac34;\frac54;\frac MN\right)
  =A_2\omega_d^2
\end{split}
\end{equation*}
hold for
$$ \small\extrarowheight3pt
\begin{array}{r|rrcc} \hline\hline
  d & M & N & A_1 & A_2 \\ \hline
-120 & -7^4 & 3^3\cdot5^3 & \displaystyle\frac12
  \sqrt[8]{45}\sqrt{12+2\sqrt{30}} & \displaystyle\frac{45}7 \\
 -52 & 2^2\cdot3^7 & 5^6 & \displaystyle\frac12
  \sqrt[4]5\sqrt{8+2\sqrt{13}} & \displaystyle\frac{25}6 \\
-132 & 2^4\cdot11^2 & 5^6 & \displaystyle
  \frac12\sqrt[8]{75}\sqrt{12+2\sqrt{33}}
     & \displaystyle\frac{75}{2\sqrt{22}} \\
 -43 & -3^7\cdot7^4 & 2^{10}\cdot5^6 & \displaystyle
  \frac12\sqrt[4]{10}\sqrt{7+\sqrt{43}} & \displaystyle
  \frac{100}{21} \\
 -88 & 3^7\cdot7^4 & 5^6\cdot11^3 & \displaystyle
  \frac12\sqrt[8]{275}\sqrt{10+2\sqrt{22}} & \displaystyle
  \frac{275}{21\sqrt2} \\
-312 & 7^4\cdot23^4 & 5^6\cdot11^6 & \displaystyle\frac12
  \sqrt[8]3\sqrt[4]{55}\sqrt{18+2\sqrt{78}}  
  & \displaystyle\frac{9075}{161\sqrt2} \\
-148 & 2^2\cdot3^7\cdot7^4\cdot11^4 & 5^6\cdot17^6
     & \displaystyle\frac12\sqrt[4]{85}\sqrt{14+2\sqrt{37}}
     & \displaystyle \frac{7225}{231}  \\
-232 & -3^7\cdot7^4\cdot11^4\cdot19^4 & 5^6\cdot23^6\cdot29^3
     & \displaystyle\frac12\sqrt[8]{29}\sqrt[4]{115}\sqrt{16+2\sqrt{58}}
     & \displaystyle\frac{383525}{4389} \\
-708 & 2^8\cdot7^4\cdot11^4\cdot47^4\cdot59^2 & 5^6\cdot17^6\cdot29^6
     & \displaystyle\frac12\sqrt[8]3\sqrt[4]{2465}\sqrt{30+2\sqrt{177}}
     & \displaystyle\frac{18228675}{3619\sqrt{118}} \\
-163 & -3^{11}\cdot7^4\cdot19^4\cdot23^4
     &  2^{10}\cdot5^6\cdot11^6\cdot17^6
     & \displaystyle\frac12\sqrt[4]{1870}\sqrt{13+\sqrt{163}}
     & \displaystyle\frac{3496900}{27531} \\
\hline\hline
\end{array}
$$
Also,
\begin{equation*}
\begin{split}
  &{}_2F_1\left(\frac1{24},\frac7{24};\frac56;\frac MN\right)
  =B_1\frac{\omega_d}{\omega_{-3}}, \\
  &{}_3F_2\left(\frac14,\frac12,\frac34;\frac56,\frac76;\frac MN\right)
  =B_2\omega_d^2
\end{split}
\end{equation*}
hold for
$$ \small\extrarowheight5pt
\begin{array}{r|rrcc} \hline\hline
 d & M & N & B_1 & B_2  \\ \hline
-84 & 3^3 & 2^2\cdot7^2 & \sqrt[12]{56}\sqrt{3+\sqrt7} &
2\sqrt{\displaystyle\frac{14}3} \\
-40 & -5^3 & 3^7 & \displaystyle\sqrt[12]{\frac43}\sqrt{2\sqrt3+\sqrt{10}} & \displaystyle\frac6{\sqrt5} \\
-51 & 2^{10} & 7^4 & \displaystyle\frac12\sqrt[6]7\sqrt{10+2\sqrt{17}}
  & \displaystyle \frac72 \\
-19 & -2^{10} & 3^7 & \displaystyle\frac12\sqrt[12]{\frac13}
  \sqrt{6\sqrt3+2\sqrt{19}}
  & \displaystyle \frac32 \\
-168 & 5^6 & 7^2\cdot11^4 & \sqrt[12]7\sqrt[6]{22}\sqrt{4+\sqrt{14}}
    & \displaystyle\frac{22}5\sqrt7 \\
-228 & -3^6\cdot5^6 & 2^6\cdot7^4\cdot19^2
     & \sqrt[12]{38}\sqrt[6]{28}\sqrt{5+\sqrt{19}}
  & \displaystyle\frac{28}{15}\sqrt{114} \\
-123 & 2^{10}\cdot5^6 & 7^4\cdot19^4 & \displaystyle\frac12
  \sqrt[6]{133}\sqrt{14+2\sqrt{41}}
     & \displaystyle \frac{133}{10} \\
-67 & -2^{16}\cdot5^6 & 3^7\cdot7^4\cdot11^4
& \displaystyle\sqrt[12]{\frac13}\sqrt[6]{\frac{77}{8}}\sqrt{5\sqrt3+\sqrt{67}}
  & \displaystyle\frac{231}{20} \\
-372 & 3^3\cdot5^6\cdot11^6 & 2^2\cdot7^4\cdot19^4\cdot31^2
     & \sqrt[12]{62}\sqrt[6]{266}\sqrt{7+\sqrt{31}}
  & \displaystyle\frac{266}{55}\sqrt{186}\\
-408 & 3^6\cdot5^6\cdot17^3 & 7^4\cdot11^4\cdot31^4
     & \sqrt[6]{4774}\sqrt{6+\sqrt{34}} & \displaystyle\frac{4774}{15\sqrt{17}} \\
-267 & 2^{16}\cdot5^6\cdot11^6 & 7^4\cdot31^4\cdot43^4
     & \displaystyle\frac12\sqrt[6]{9331}\sqrt{22+2\sqrt{89}}
     & \displaystyle\frac{9331}{110} \\
\hline\hline
\end{array}
$$
\end{Theorem}


\begin{Remark} Let $F_1(s)={}_2F_1(1/24,5/24;3/4;s)$,
  $G_1(t)={}_2F_1(1/24,7/24;5/6;t)$, and
\begin{equation*}
\begin{split}
  F_2(s)&={}_2F_1(7/24,11/24;5/4;s)={}_3F_2(1/3,1/2,2/3;3/4;5/4;s)/F_1(s) \\
  G_2(t)&={}_2F_1(5/24,11/24;7/6;t)={}_3F_2(1/4,1/2,3/4;5/6,7/6;t)/G_1(t)
\end{split}
\end{equation*}
be the hypergeometric functions in Theorem \ref{theorem: main theorem}.
The Ramanujan-type identities obtained in
\cite{Yang-Ramanujan} can be written as
\begin{equation*}
\begin{split}
  \left(R_1s\frac{d}{ds}F_1(s)^2+R_2F_1(s)^2\right)\Big|_{s=M/N}
&=\sqrt{R_3}|M|^{3/4}N^{1/4}C_1, \\
  \left(R_1s\frac d{ds}F_2(s)^2+(R_1/2+R_2)F_2(s)\right)\Big|_{s=M/N}
&=\sqrt{R_3}|M|^{1/4}N^{3/4}C_1^{-1},
\end{split}
\end{equation*}
and
\begin{equation*}
\begin{split}
  \left(R_1t\frac d{dt}G_1(t)^2+R_2G_1(t)^2\right)\Big|_{t=M/N}
&=\sqrt{R_3}|M|^{2/3}N^{1/3}C_2, \\
  \left(R_1t\frac d{ds}G_2(t)^2+(R_1/3+R_2)G_2(t)\right)\Big|_{t=M/N}
&=\sqrt{R_3}|M|^{1/3}N^{2/3}C_2^{-1},
\end{split}
\end{equation*}
for some rational numbers $R_1,R_2,R_3$ depending on $d$, where
$$
  C_1=\frac{4}{\sqrt[4]{12}}\frac\pi{\Omega_{-4}^{2}}
     =\frac4{\sqrt[4]{12}}\frac{\Gamma(3/4)^2}{\Gamma(1/4)^2}, \qquad
  C_2=\frac{3}{\sqrt[6]2}\frac\pi{\Omega_{-3}^{2}}
     =\frac3{\sqrt[6]2}\frac{\Gamma(2/3)^3}{\Gamma(1/3)^3}.
$$
Combining these identities with the formulas in Theorem \ref{theorem:
  main theorem}, we obtain special values for the functions
\begin{equation*}
\begin{split}
  \frac d{ds}F_1(s)^2&=\frac d{ds}{}_3F_2\left(
  \frac1{12},\frac14,\frac5{12};\frac12,\frac34;s\right)
 =\frac5{216}{}_3F_2\left(\frac{13}{12},\frac54,\frac{17}{12};\frac32,\frac74;
  s\right) \\
  \frac d{ds}F_2(s)^2&=\frac d{ds}{}_3F_2\left(
  \frac7{12},\frac34,\frac{11}{12};\frac32,\frac54;s\right)
 =\frac{77}{360}{}_3F_2\left(\frac{19}{12},\frac74,\frac{23}{12};
  \frac52,\frac94;s\right),
\end{split}
\end{equation*}
For instance, for $d=-120$, we have
\begin{equation*}
\begin{split}
  _3F_2\left(\frac{13}{12},\frac54,\frac{17}{12};\frac32,\frac74;
  -\frac{7^4}{15^3}\right)&=\frac{3^6\cdot5^{9/4}}
  {2\cdot7^3\cdot19\cdot\omega_{-4}^2}
  \left((4\sqrt3+2\sqrt{10})\omega_{-120}^2-\sqrt{\frac32}\right) \\
  _3F_2\left(\frac{19}{12},\frac74,\frac{23}{12};\frac52,
  \frac94;-\frac{7^4}{15^3}\right)
  &=\frac{3^7\cdot5^{23/4}\cdot\omega_{-4}^2}{7^7\cdot11\cdot19}
   \left(242(2\sqrt3-\sqrt{10})\omega_{-120}^2-7\sqrt{\frac32}\right).
\end{split}
\end{equation*}
There are similar formulas for the functions
$$
  _3F_2\left(\frac{13}{12},\frac43,\frac{19}{12};\frac{11}6,\frac53;t
  \right), \quad
  _3F_2\left(\frac{17}{12},\frac53,\frac{23}{12};\frac{13}6,\frac73;t
  \right),
$$
such as
\begin{equation*}
\begin{split}
  _3F_2\left(\frac{13}{12},\frac43,\frac{19}{12};\frac{11}6,\frac53;
  \frac{27}{196}\right)
  &=\frac{2^4\cdot5\cdot7^{7/6}}{3^2\cdot13\cdot\omega_{-3}^2}\left(\frac4{\sqrt3}
   -\sqrt2(3+\sqrt7)\omega_{-84}^2\right), \\
  _3F_2\left(\frac{17}{12},\frac53,\frac{23}{12};\frac{13}6;\frac73;
  \frac{27}{196}\right)
  &=\frac{2^5\cdot7^{23/6}\cdot\omega_{-3}^2}{3^4\cdot5\cdot11\cdot13}
    \left(4\sqrt3-55\sqrt2(3-\sqrt7)\omega_{-84}^2\right).
\end{split}
\end{equation*}
\end{Remark}

\begin{Remark}
Notice that the numbers $A_1$ in the first table are all of the form
$A^{1/8}(a+\sqrt{|d|})^{1/2}$ for some positive integer $a$ and some
rational number $A$ whose denominator is $2$ or $4$. In other words,
the special values $_2F_1(1/24,5/24;3/4;M/N)$ possess a certain
integrality property. This integrality property is a consequence of
Schofer's work \cite{Schofer} and our explicit realization of modular
forms as Borcherds form. On the other hand, if we can somehow manage
to prove this integrality property without using Borcherds forms, then
to obtain the identities in Theorem \ref{theorem: main theorem}, we can
just evaluate the hypergeometric functions to a high precision and identify
the integers. Note that the prime factors of the numerator of $A$ are
either $2$ or prime factors of $N$. This suggests that it may be
possible to prove the integrality property using the moduli
interpretation of the Shimura curve $X_0^6(1)$.
\end{Remark}

\begin{Remark} Note that the proof of Theorem \ref{theorem: 0} is
  certainly valid for other Shimura curves $X_0^D(N)/W$, $W$ being a
  subgroup of the Atkin-Lehner groups, or even
  Shimura curves over totally real fields. However, other than the
  cases of arithmetic triangle groups, as classified in
  \cite{Takeuchi}, there are only a very limited number of Shimura
  curves whose Schwarzian differential equations are known (see
  \cite{Elkies,Tu-Schwarzian}).

  To obtain analogues of Theorem \ref{theorem: main theorem} for
  $X_0^D(N)/W$, one will need a method to construct Borcherds forms
  systematically. This is recently addressed in \cite{Guo-Yang}, so
  there is no problem in evaluating modular forms on
  $X_0^D(N)/W$ at CM-points. However, we remark that this only
  translates to analogues of the $_2F_1$-evaluations. To obtain
  analogues of the $_3F_2$-evaluations, we will need to determine the
  constant $C$ such that the linear combination $f_1+Cf_2$ of two
  solutions $f_1$ and $f_2$ of the Schwarzian differential equation is
  a modular form. In general, this is a difficult problem.
  (For the case of $X_0^6(1)/W_6$, the constant $C$ is determined by
  using Gauss' formula
  $_2F_1(a,b;c;1)=\Gamma(c)\Gamma(c-a-b)/\Gamma(c-a)\Gamma(c-b)$.)

  If one wishes to further generalize Theorem \ref{theorem: main
    theorem} to Shimura curves over totally real fields, one will need
  the theory of Borcherds forms over totally real fields, developed
  recently by Bruinier and Yang \cite{Bruinier-Yang,Bruinier}.
  As far as we can see, it should in principle be possible to obtain
  explicit evaluations at least for the case of arithmetic triangle groups. We
  leave this problem for future investigation.
\end{Remark}

\begin{Remark}
Notice that if a prime $p$ divides $M$, then the hypergeometric series
appearing in Theorem \ref{theorem: main theorem} converges
$p$-adically and one may wonder what the limit is. Our computation
suggests the following $p$-adic evaluation.

For a prime $p$, let $\Gamma_p(x)$ be the $p$-adic Gamma function
defined by
$$
  \Gamma_p(n)=(-1)^n\prod_{0<j<n,p\nmid j}j
$$
for positive integers $n$ and extended continuously to $\Z_p$, and for
a fundamental discriminant $d<0$, set
$$
  \omega_{d,p}=\prod_{a=1}^{|d|-1}\Gamma_p\left(\frac{a}{|d|}
  \right)^{\chi_d(a)\mu_d/8h_d}.
$$
Consider the two hypergeometric functions in the first set of
identities in Theorem \ref{theorem: main theorem}. Other than the
cases $d=-52$ and $d=-132$, the series converge $7$-adically. 
Then the numerical data suggest that
\begin{equation*}
\begin{split}
 &{}_2F_1\left(\frac1{24},\frac5{24};\frac34;\frac MN\right)
 =A_1\frac{\omega_{d,7}}{\omega_{-4,7}}, \\
 &{}_3F_2\left(\frac13,\frac12,\frac23;\frac34,\frac54;\frac MN\right)
 =A_2\omega_{d,7}^2
\end{split}
\end{equation*}
hold with the same $M$ and $N$, and
$$ \extrarowheight3pt
\begin{array}{c|cc} \hline\hline
d & A_1 & A_2 \\ \hline
-120 & \displaystyle\frac12\sqrt[8]{-\frac1{125}}
  \sqrt{\sqrt{-3}+\sqrt{-10}} & \displaystyle\frac3{8\sqrt2} \\
-43 & \displaystyle\sqrt[4]{-10}\sqrt{\frac{1+\sqrt{43}}{43}}
  & -\displaystyle\frac{100}{129} \\
-88 & \displaystyle\frac12\sqrt[8]{-\frac1{11}}\sqrt[4]{\frac5{11}}
  \sqrt{1+\sqrt{22}} & \displaystyle\frac{25}{24\sqrt2} \\
-312 & \displaystyle\sqrt[8]{-\frac13}\sqrt[4]{-\frac{55}8}
   \sqrt{\frac{\sqrt{-6}+\sqrt{-13}}{13}}
  & \displaystyle-\frac{3025}{2392} \\
-148 & \displaystyle\sqrt[4]{85}\sqrt{\frac{4+\sqrt{37}}{74}}
  & \displaystyle -\frac{7225}{9768} \\
-232 & \displaystyle\sqrt[8]{\frac1{29}}\sqrt[4]{-\frac{115}{232}}
  \sqrt{2\sqrt2+\sqrt{29}} & \displaystyle \frac{13225}{5016} \\
-708 & \displaystyle \sqrt[8]{-\frac13}
   \sqrt[4]{2465}\sqrt{\frac{\sqrt{-3}+\sqrt{-59}}{118}}
  & \displaystyle\frac{6076225}{244024\sqrt{-59}} \\
-163 & \displaystyle\sqrt[4]{1870}
   \sqrt{\frac{11+\sqrt{163}}{163}}
  & -\displaystyle\frac{3496900}{641079} \\
\hline\hline
\end{array}
$$
(There are many places where we need to take square roots of $p$-adic
numbers. The table above means that after taking suitable choices of
square roots, the identities hold conjecturally.)

Note that for a prime $p$ and a fundamental discriminant $d<0$,
the $p$-adic number $\omega_{d,p}^2$
appears in the matrix representation of the Frobenius automorphism on
the de Rham cohomology $H_{\mathrm{dR}}^1(E/\overline\Q)\otimes
K_{\mathfrak p}$ for an elliptic curve $E$ over $\overline\Q$ with CM
by $\Q(\sqrt d)$, where $\mathfrak p$ is the prime of $\overline\Q$
lying over $p$ and $K_{\mathfrak p}$ is the algebraic closure of
$\Q_p$ in the completion of $\overline\Q$ at $\mathfrak p$.
(See \cite[Theorem 3.15]{Ogus}.) Note also that if the prime $p$ splits in
$\Q(\sqrt d)$, then $w_{d,p}$ is algebraic over $\Q$ since a suitable
power of $\omega_{d,p}$ appears as the value of a certain $p$-adic
Gaussian sum. (See \cite[Theorem 1.12]{Gross-Koblitz}.) On the other
hand, it is expected that when $p$ is inert in $\Q(\sqrt d)$,
$\omega_{d,p}$ is transcendental over $\Q$. In our conjectural
$7$-adic formulas mentioned above, since the prime $7$ is always 
inert in $\Q(\sqrt d)$ (which is a consequence of Theorem 3.6
of \cite{Schofer}), we expect that $\omega_{d,7}$ is transcendental
over $\Q$ for $d$ given in the list.
\end{Remark}
\end{section}

\begin{section}{Realization of modular forms in terms of Schwarzian
    differential equations}
\label{section: Schwarzian}

Here we briefly explain the realization of modular forms on Shimura
curves using solutions of Schwarzian differential equations. For
details, see \cite{Yang-Schwarzian}.

Assume that a Shimura curve $X$ has genus $0$ with elliptic points and
cusps $\tau_1,\ldots,\tau_r$ of order $e_1,\ldots,e_r$, respectively.
(Here we set $e_j=\infty$ if $\tau_j$ is a cusp.) Let $t(\tau)$ be a
Hauptmodul for $X$ and set $a_j=t(\tau_j)$. Then Theorem 4 of
\cite{Yang-Schwarzian} shows that a basis for the space of modular
forms of even weight $k\ge 4$ is
\begin{equation} \label{equation: basis}
  t'(\tau)^{k/2}t(\tau)^j\prod_{i=1,a_i\neq\infty}^r
 (t(\tau)-a_i)^{-\lfloor k(1-1/e_i)/2\rfloor}, \quad j=0,\ldots, d_k-1,
\end{equation}
where
$$
  d_k=1-k+\sum_{j=1}^r\left\lfloor\frac k2\left(1-\frac1{e_j}\right)
  \right\rfloor
$$
is the dimension of the space of modular forms of weight $k$ on $X$.

Now it is easy to check that $t'(\tau)$ is a meromorphic modular form
of weight $2$ on $X$. Thus, $t'(\tau)^{1/2}$ and $\tau
t'(\tau)^{1/2}$, as functions of $t$, are solutions of a certain
second-order linear differential equation with rational functions in
$t$ as coefficients. (See \cite[Theorem 5.1]{Stiller} or
\cite[Theorem 1]{Yang}. Here the coefficients of the differential
equation are rational functions because $t$ is a Hauptmodul.) In
fact, this differential equation is
\begin{equation} \label{equation: Schwarzian}
  \frac{d^2}{dt^2}F+Q(t)F=0,
\end{equation}
where
$$
  Q(t)=-\frac12\frac{\{t,\tau\}}{t'(\tau)^2}, \qquad
  \{t,\tau\}=\frac{t'''(\tau)}{t'(\tau)}-\frac32\left(
  \frac{t''(\tau)}{t'(\tau)}\right)^2.
$$
Because $\{t,\tau\}$ is classically known as the Schwarzian
derivative, we call the differential equation satisfied by
$t'(\tau)^{1/2}$ and $t(\tau)$ the \emph{Schwarzian differential
  equation} associated to the Shimura curve. If we let $\{f_1,f_2\}$
be a basis for the solution of \eqref{equation: Schwarzian}, then we
have $t'(\tau)=(c_1f_1+c_2f_2)^2$ for some complex numbers $c_1$ and
$c_2$. Substituting this into \eqref{equation: basis}, we obtain
realization of modular forms in terms of solutions of Schwarzian
differential equations.

When a Shimura curve is of genus zero and has precisely three elliptic points
or cusps, the Schwarzian differential equation is essentially a
hypergeometric differential equation. In particular, for the curve
$X=X_0^6(1)/W_6$, we can realize modular forms on $X$ in terms of
hypergeometric functions as follows.

We let $B=\Q+\Q I+\Q J+\Q IJ$ with $I^2=-1$, $J^2=3$, and $IJ=-JI$, be
the quaternion algebra of discriminant $6$ over $\Q$ and choose the
embedding $\iota: B\hookrightarrow M(2,\R)$ to be the one defined by
$$
  \iota(I)=\M0{-1}10, \qquad \iota(J)=\M{\sqrt3}00{-\sqrt3}.
$$
Fix a maximal order $\O=\Z+\Z I+\Z J+\Z(1+I+J+IJ)/2$ in $B$ and choose
representatives of CM-points of discriminants $-3$, $-4$, and $-24$ to be
$$
  P_{-3}=\frac{-1+i}{1+\sqrt3}, \quad P_{-4}=i, \quad
  P_{-24}=\frac{(\sqrt6-\sqrt2)i}2.
$$
They are the elliptic points of orders $6$, $4$, and $2$,
respectively. A fundamental domain is given by

\centerline{\epsfig{file=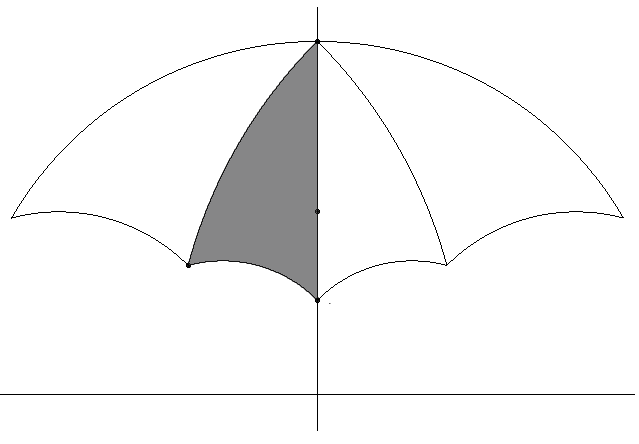,height=2.1in,width=3in}}

\noindent Here the grey area represents a fundamental domain for
$X_0^6(1)/W_6$. The four marked points on the boundary are
$P_{-4}$, $P_{-3}$, $P_{-24}$, and $(2-\sqrt3)i$.

We have the following bases for the spaces of
modular forms on $X_0^6(1)/W_6$.

\begin{Proposition} \label{proposition: realization 1}
Let $s$ be the Hauptmodul on $X=X_0^6(1)/W_6$ determined by
$s(P_{-4})=0$, $s(P_{-24})=1$, and $s(P_{-3})=\infty$. Then for an
even integer $k\ge 4$, a basis for the space $S_k(X)$ of modular forms
of weight $k$ on $X$ is
$$
  s^{\{3k/8\}}(1-s)^{\{k/4\}}s^j
  \left({}_2F_1\left(\frac1{24},\frac5{24};\frac34;s\right)
  +\frac1{\sqrt[4]{12}\omega_{-4}^2}s^{1/4}
  {}_2F_1\left(\frac7{24},\frac{11}{24};\frac54;s\right)\right)^k,
$$
$j=0,\ldots,d_k-1$, where $d_k=\dim
S_k(X)=1-k+\FL{k/4}+\FL{3k/8}+\FL{5k/12}$.
 
Also, let $t=1/s$. Then a basis for $S_k(X)$ is
$$
  t^{\{5k/12\}}(1-t)^{\{k/4\}}t^j
  \left({}_2F_1\left(\frac1{24},\frac7{24};\frac56;t\right)
  -\frac{e^{-2\pi i/8}}{\sqrt[6]2\omega_{-3}^2}t^{1/6}
  {}_2F_1\left(\frac5{24},\frac{11}{24};\frac76;t\right)\right)^k,
$$
$j=0,\ldots,d_k-1$.
\end{Proposition}

\begin{proof} The first part is the content of Lemmas 3 and 4 of
  \cite{Yang-Ramanujan}. For the second part, the proof of Lemma 14 of
  \cite{Yang-Schwarzian} shows that
  \begin{equation} \label{equation: t'}
    t'(\tau)=\frac{6t^{5/6}(1-t)^{1/2}}{C(P_{-3}-\overline P_{-3})}
    \left({}_2F_1\left(\frac1{24},\frac7{24};\frac56;t\right)
   -Ct^{1/6}{}_2F_1\left(\frac5{24},\frac{11}{24};\frac76;t\right)\right)^2,
  \end{equation}
  where
  $$
    C=\frac{P_{-24}-P_{-3}}{P_{-24}-\overline P_{-3}}
      \frac{\Gamma(5/6)\Gamma(17/24)\Gamma(23/24)}
      {\Gamma(7/6)\Gamma(13/24)\Gamma(19/24)}.
  $$
  Now
  $$
    \frac{P_{-24}-P_{-3}}{P_{-24}-\overline P_{-3}}
   =(1-i)\left(1-\frac1{\sqrt2}\right)=e^{-2\pi i/8}(\sqrt2-1).
  $$
  Also, by Euler's reflection formula and
  Gauss's multiplication formula, we have
  \begin{equation*}
  \begin{split}
    \left(\frac{\Gamma(17/24)\Gamma(23/24)}
    {\Gamma(13/24)\Gamma(19/24)}\right)^2
  &=\frac{\Gamma(17/24)\Gamma(23/24)\Gamma(5/24)\Gamma(11/24)}
    {\Gamma(13/24)\Gamma(19/24)\Gamma(1/24)\Gamma(7/24)} \\
  &\qquad\qquad\times
    \frac{\sin(5\pi/24)\sin(11\pi/24)}{\sin(\pi/24)\sin(7\pi/24)} \\
  &=4^{-2/3}\frac{\Gamma(5/6)}{\Gamma(1/6)}(3+2\sqrt2)
  \end{split}
  \end{equation*}
  and
  $$
    \Gamma\left(\frac13\right)\Gamma\left(\frac56\right)
   =(2\pi)^{1/2}2^{-1/6}\Gamma\left(\frac23\right).
  $$
  From these, we deduce that
  \begin{equation*}
  \begin{split}
    \frac{\Gamma(5/6)}{\Gamma(7/6)}
    \frac{\Gamma(17/24)\Gamma(23/24)}{\Gamma(13/24)\Gamma(19/24)}
  &=6\cdot2^{-2/3}(\sqrt2+1)\frac{\Gamma(5/6)^{3/2}}{\Gamma(1/6)^{3/2}} \\
  &=6\cdot2^{-2/3}(\sqrt2+1)\frac1{(2\pi)^{3/2}}\Gamma\left(\frac56\right)^3\\
  &=6\cdot2^{-7/6}(\sqrt2+1)\frac{\Gamma(2/3)^3}{\Gamma(1/3)^3}
   =\frac{\sqrt2+1}{\sqrt[6]2\omega_{-3}^2}
  \end{split}
  \end{equation*}
  and hence
  $$
    C=\frac{e^{-2\pi i/8}}{\sqrt[6]2\omega_{-3}^2}.
  $$
  Then from \eqref{equation: basis}, we conclude that the second set
  of functions in the lemma forms a basis for $S_k(X)$.
\end{proof}

For general Shimura curves, we can determine
Schwarzian differential equations using Propositions 5 and 6 of
\cite{Yang-Schwarzian} and explicit covers of Shimura curves. In
\cite{Tu-Schwarzian}, Tu determines Schwarzian differential equations
for the cases when $X_0^D(1)/W_D$ and $X_0^D(N)/W_D$ both have genus
zero.

We now give a proof of Theorem \ref{theorem: 0}.

\begin{proof}[Proof of Theorem \ref{theorem: 0}]
  Here we only prove the second half of the theorem; the proof of the
  first half is similar and is omitted.

  Since $t(\tau)$ is a Hauptmodul that takes rational values at
  three distinct CM-points, it takes algebraic values at all
  CM-points. Thus, by \cite[Theorem 7.1]{Shimura-periods} and
  \cite[Theorem 1.2 and (1.4) of Chapter 3]{Yoshida}, the value of
  $t'(\tau)$ at a CM-point of discriminant $d$ is an algebraic
  multiple of $\omega_{d_0}^2$. Then, from \eqref{equation: t'}, we
  see that
  $$
    {}_2F_1\left(\frac1{24},\frac7{24};\frac56;t(\tau_d)\right)
   -\frac{e^{-2\pi i/8}}{\sqrt[6]2\omega_{-3}^2}t(\tau_d)^{1/6}
    {}_2F_1\left(\frac5{24},\frac{11}{24};\frac76;t(\tau_d)\right)
    \in\frac{\omega_{d_0}}{\omega_{-3}}\cdot\overline\Q.
  $$
  Without loss of generality, we may assume that $\tau_d$ lies in the
  fundamental domain depicted earlier. Then Equation (22) of
  \cite{Yang-Schwarzian} implies that
  $$
    \frac{{}_2F_1(5/24,11/24;7/6;t(\tau_d))}
    {{}_2F_1(1/24,7/24;5/6;t(\tau_d))}\in\omega_{-3}^2\cdot\overline\Q.
  $$
  It follows that
  $$
    {}_2F_1\left(\frac1{24},\frac7{24};\frac56;t(\tau_d)\right)
    \in\frac{\omega_{d_0}}{\omega_{-3}}\cdot\overline\Q
  $$
  and
  \begin{equation*}
  \begin{split}
    {}_3F_2\left(\frac14,\frac12,\frac34;\frac56,\frac76;
    t(\tau_d)\right)
  &={}_2F_1\left(\frac1{24},\frac7{24};\frac56;t(\tau_d)\right)
    {}_2F_1\left(\frac5{24},\frac{11}{24};\frac76;t(\tau_d)\right) \\
  &\in\omega_{d_0}^2\cdot\overline\Q.
  \end{split}
  \end{equation*}
  This proves the theorem.
\end{proof}
\end{section}

\begin{section}{Realization of modular forms as Borcherds forms}
\label{section: Borcherds forms}

We first give a quick introduction to Borcherds forms. For details,
see \cite{Borcherds,Borcherds-Duke}.

Let $L$ be an even lattice with a symmetric bilinear form
$\gen{\cdot,\cdot}$ of signature $(b^+,b^-)$,
$L^\vee=\{\gamma\in L\otimes\Q: \gen{\gamma,\eta}\in\Z\text{ for all
}\eta\in L\}$ its dual lattice, and $\{e_\eta:\eta\in L^\vee/L\}$ the
standard basis for the vector space $\C[L^\vee/L]$. Let
$$
  \wt\SL(2,\Z)=\left\{\left(\M abcd,\pm\sqrt{c\tau+d}\right):
  \M abcd\in\SL(2,\Z)\right\}
$$
be the metaplectic double cover of $\SL(2,\Z)$, which is generated by
$$
  S=\left(\M0{-1}10,\sqrt\tau\right), \qquad
  T=\left(\M1101,1\right).
$$
Associated to the lattice $L$, we have the \emph{Weil representation}
$\rho_L:\wt\SL(2,\Z)\to\GL(\C[L^\vee/L])$ defined by
\begin{equation*}
\begin{split}
  \rho_L(T) e_\eta&=e^{-2\pi i\gen{\eta,\eta}/2} e_\eta, \\
  \rho_L(S) e_\eta&=\frac{e^{2\pi i(b^+-b^-)/8}}{\sqrt{|L^\vee/L|}}
  \sum_{\gamma\in L^\vee/L}e^{2\pi i\gen{\eta,\gamma}}e_\gamma.
\end{split}
\end{equation*}
A holomorphic function $F:\H^+\to\C[L^\vee/L]$ is said to be a
\emph{weakly holomorphic vector-valued modular form} of weight
$k\in\frac12\Z$ and type $\rho_L$ if it satisfies
$$
  F\left(\frac{a\tau+b}{c\tau+d}\right)
 =(c\tau+d)^k\rho_L\left(\M abcd,\sqrt{c\tau+d}\right)F(\tau)
$$
for all $\tau\in\H^+$ and all $\SM abcd\in\SL(2,\Z)$ and the principal
part of its Fourier expansion
$F(\tau)=\sum_\eta(\sum_{m\in\Q}c_\eta(m)q^m)e_\eta$, $q=e^{2\pi
  i\tau}$, has finitely many terms, i.e., the number of pairs
$(\eta,m)$ with $m<0$ and $c_\eta(m)\ne 0$ is finite.

For $k=\Q$, $\R$, or $\C$, let $V(k)=L\otimes k$ and extend the
definition of $\gen{\cdot,\cdot}$ to $V(k)$ by linearity. Define
the orthogonal groups
$$
  O_V(\R)=\{\sigma\in\GL(V(\R)):\gen{\sigma x,\sigma y}=\gen{x,y}
  \text{ for all }x,y\in V(\R)\}
$$
and
$$
  O_V^+(\R)=\{\sigma\in O_V(\R):\operatorname{spin}\sigma=\sgn\det\sigma\},
$$
where if $\sigma$ is equal to the product of $n$ reflections with
respect to the vectors $v_1,\ldots,v_n$, then its spinor norm is
defined by $\operatorname{spin}\sigma=(-1)^n\prod_{i=1}^n\sgn\gen{v_i,v_i}$.
Let also
$$
  O_L=\{\sigma\in O_V(\R):~\sigma(L)=L\}, \qquad
  O_L^+=O_L\cap O_V^+(\R).
$$
(Note that the definition of spinor norms is different from that of
\cite{Borcherds} since the bilinear form in our setting differs from
that of \cite{Borcherds} by a factor of $-1$.) 

From now on, we assume that the signature of $L$ is $(b,2)$. Let
$\Gr(V(\R))$ be the Grassmanian of oriented negative $2$-planes in
$V(\R)$. For an element $A$ in $\Gr(V(\R))$, we can find an oriented
basis $\{x,y\}$ for $A$ with $\gen{x,x}=\gen{y,y}=-1$ and
$\gen{x,y}=0$. Let $z=x+iy\in V(\C)$. Then we have $\gen{z,z}=0$ and
$\gen{z,\overline z}<0$. In fact, it is easy to show that $\Gr(V(\R))$
can be identified with the set
$$
  K=\{z\in V(\C):\gen{z,z}=0,\gen{z,\overline z}<0\}/\C^\times.
$$
The set $K$ has two connected components, which amount to the two
choices of continuously varying orientation of negative $2$-planes in
$V(\R)$. Pick one of them to be $K^+$. Then the orthogonal group
$O_V^+(\R)$ acts transitively on $K^+$. Let
$$
  \widetilde K^+=\{z\in V(\C):\gen{z,z}=0,\gen{z,\overline z}<0,
  [z]\in K^+\}.
$$
Then for a subgroup $\Gamma$ of $O_L^+$, a meromorphic function
$\Psi:\widetilde K^+\to\P^1(\C)$ is called a modular form of weight
$k$ with character $\chi$ on $\Gamma$ if $\Psi$ satisfies
\begin{enumerate}
\item $\Psi(cz)=c^{-k}\Psi(z)$ for all $c\in\C^\times$ and
  $z\in\widetilde K^+$, and
\item $\Psi(gz)=\chi(g)\Psi(z)$ for all $g\in\Gamma$ and
  $z\in\widetilde K^+$.
\end{enumerate}

\begin{theorem}[{\cite[Theorem 13.3]{Borcherds}}] Let $L$ be an even
  lattice of signature $(b,2)$ and $F(\tau)$ be a weakly holomorphic
  vector-valued modular form of weight $1-b/2$ and type $\rho_L$ with
  Fourier expansion $F(\tau)=\sum_{\eta\in
    L^\vee/L}F_\eta(\tau)e_\eta=\sum_{\eta}(\sum_{m\in\Q}
  c_\eta(m)q^m)e_\eta$. Suppose that $c_\eta(m)\in\Z$ whenever $m\le
  0$. Then there corresponds a meromorphic function $\Psi_F(z)$,
  $z\in\widetilde K^+$, with the following properties.
  \begin{enumerate}
  \item $\Psi_F(z)$ is a meromorphic modular form of weight $c_0(0)/2$
    on the group
  $$
    O_{L,F}^+=\{\sigma\in O_L^+:~F_{\sigma\eta}=F_\eta\text{ for all }
    \eta\in L^\vee/L\}
  $$
  with respect to some unitary character $\chi$ of $O_{L,F}^+$.
  \item The only zeros or poles of $\Psi_F(z)$ lie on the rational
    quadratic divisor $\lambda^\perp=\{z\in\widetilde
    K^+:\gen{z,\lambda}=0\}$ for $\lambda\in L$,
    $\gen{\lambda,\lambda}>0$ and are of order
    $$
      \sum_{0<r\in\Q,r\lambda\in L^\vee}c_{r\lambda}(-r^2\gen{\lambda,\lambda}/2).
    $$
  \end{enumerate}
\end{theorem}

We call the function $\Psi_F(z)$ the \emph{Borcherds form} associated
to $F(\tau)$.

We now explain the idea of realizing modular forms on Shimura curves
in terms of Borcherds forms. Even though this idea has been used in
\cite{Errthum}, it seems to us that some key properties were not
explained very concretely in \cite{Errthum}. For instance, it was not
explained in \cite{Errthum} why the characters associated to the
Borcherds forms constructed therein are trivial. Therefore,
it is worthwhile to explain this approach in some
details.

Let $\O$ be an Eichler order of level $N$ in an indefinite quaternion
algebra $B$ of discriminant $D$ over $\Q$, $(N,D)=1$, $\O_1$ be
the group of norm-one elements in $\O$, and
$$
  L=\{\alpha\in\O:\tr(\alpha)=0\}
$$
be the set of elements of trace zero in $\O$, where $\tr(\alpha)$ and
$\n(\alpha)$ denote the trace and the norm of $\alpha$, respectively.
By setting $\gen{\alpha,\beta}=\tr(\alpha\beta')$, $L$ becomes a
lattice of signature $(1,2)$, where $\beta'$ denote the quaternionic
conjugate of $\beta$ in $B$. We now determine $O_L$ and $O_L^+$.

By the Cartan-Dieudonn\'e theorem, every isometry $\sigma$ in
$O_V(\R)$ is equal to the product of at most three reflections. Now
it is clear that for an element of nonzero norm $\alpha$ in $V(\R)$,
the function $\tau_\alpha:\gamma\to-\alpha\gamma\alpha^{-1}$ sends
$\alpha$ to $-\alpha$ and leaves any element of $V(\R)$ orthogonal to
$\alpha$ fixed. (Here we regard $V(\R)$ as the set of trace-zero
elements in the quaternion algebra $B\otimes\R$ and define
multiplication and inverse accordingly.) In other words, $\tau_\alpha$
is the reflection with respect to $\alpha$.
Thus, $\sigma$ has determinant $1$, i.e., $\sigma$ is the product of
an even number of reflections, if and only if $\sigma$ is the
isometry $\sigma_\beta:\gamma\to\beta\gamma\beta^{-1}$
induced by the conjugation by an element $\beta$ of nonzero norm in
$B\otimes\R$ and $\sigma$ has determine $-1$ if and only if
$\sigma=-\sigma_\beta$ for some $\beta$. From this, we deduce that
$$
  O_V(\R)=\{\sigma_\beta:\beta\in(B\otimes\R)/\R^\times,~\n(\beta)\ne 0\}
  \times\{\pm 1\}
$$
and
$$
  O_V^+(\R)=\{\sigma_\beta:\beta\in(B\otimes\R)/\R^\times,~\n(\beta)>0\}
  \times\{\pm 1\}.
$$
In addition, by the Noether-Skolem theorem, if $\sigma_\beta$,
$\beta\in B\otimes\R$, satisfies $\sigma_\beta(V(\Q))=V(\Q)$, then
$\beta$ can be chosen from $B$. It follows that
$$
  O_L=\{\sigma_\beta:\beta\in N_B(\O)/\Q^\times\}\times\{\pm 1\}
$$
and
$$
  O_L^+=\{\sigma_\beta:\beta\in N_B^+(\O)/\Q^\times\}\times\{\pm 1\},
$$
where $N_B(\O)$ denotes the normalizer of $\O$ in $B$ and $N_B^+(\O)$
is the subgroup of elements of positive norm in $N_B(\O)$.

Now assume that the quaternion algebra $B$ is represented by
$B=\JS{a,b}\Q$ with $a,b>0$. That is, $B=\Q+\Q I+\Q J+\Q IJ$ with
$I^2=a$, $J^2=b$, and $IJ=-JI$. Fix an embedding $\iota:B\to
M(2,\R)$ by
$$
  \iota: I\to\M0{\sqrt a}{\sqrt a}0, \quad
  J\to\M{\sqrt b}00{-\sqrt b}.
$$
We can show that each class in $K=\{z\in
V(\C):\gen{z,z}=0,\gen{z,\overline z}<0\}/\C^\times$ contains a
unique representative of the form
\begin{equation} \label{equation: z(tau)}
  z(\tau)=\frac{1-\tau^2}{2\sqrt a}I+\frac{\tau}{\sqrt b}J
  +\frac{1+\tau^2}{2\sqrt{ab}}IJ
\end{equation}
for some $\tau\in\H^\pm$, the union of the upper and lower
half-planes, and the mapping $\tau\to z(\tau)\mod\C^\times$ is a
bijection between $\H^\pm$ and $K$. Let $K^+$ be the image of $\H^+$
under this mapping. Now the group $N_B^+(\O)/\Q^\times$ acts on $\H^+$ 
by linear fractional transformation through the embedding $\iota$ and
also on $K^+$ by conjugation. By a straightforward computation, we can
verify that the actions are compatible. To be more concrete, for
$\alpha\in N_B^+(\O)$, if we write
$\iota(\alpha)=\SM{c_1}{c_2}{c_3}{c_4}$, then for all $\tau\in\H^+$,
we have
\begin{equation} \label{equation: compatibility}
  \alpha z(\tau)\alpha^{-1}=\frac{(c_3\tau+c_4)^2}{\n(\alpha)}
  z\left(\frac{c_1\tau+c_2}{c_3\tau+c_4}\right).
\end{equation}
Thus, if $\Psi(z)$ is a meromorphic modular form of weight $k$ on 
$O_L^+$ with character $\chi$, then the function $\psi(\tau)$ defined
by $\psi(\tau)=\Psi(z(\tau))$ is a meromorphic modular form of weight
$2k$ with character on the Shimura curve
$N_B^+(\O)\backslash\H^+$. Since the group $N_B^+(\O)/(\Q^\times\O_1)$
contains the Atkin-Lehner group, we find that $\psi(\tau)$ is a
modular form on $X_0^D(N)/W_{D,N}$, the quotient of the Shimura curve
$X_0^D(N)$ by the group $W_{D,N}$ of all Atkin-Lehner involutions.
In particular, we have the following lemma.

\begin{Lemma} \label{lemma: psi}
Let $F(\tau)=\sum_\eta(\sum_m c_\eta(m)q^m)e_\eta$ be a
weakly holomorphic vector-valued modular form of weight $1/2$ and type
$\rho_L$ such that $O_{L,F}^+=O_L^+$ and $c_\eta(m)\in\Z$ whenever
$m\le 0$. Then the function $\psi_F(\tau)$ defined by
$\psi_F(\tau)=\Psi_F(z(\tau))$ is a meromorphic modular 
form of weight $c_0(0)$ with certain unitary character $\chi$ on the
Shimura curve $X_0^D(N)/W_{D,N}$.
\end{Lemma}


We now determine the divisor of $\psi_F(\tau)$.
According to Borcherds' theorem, the divisor of $\Psi_F(z)$ is
supported on $\lambda^\perp$ for $\lambda\in L$ with positive norm
such that $c_{r\lambda}(-r^2\n(\lambda))\ne 0$ for some positive
rational number $r$. Now suppose that $\lambda$ is
such an element of $L$. The condition $\gen{\lambda,z}=0$
implies that $\lambda z\lambda^{-1}=-z=z\mod\C^\times$. That is,
$\lambda^\perp/\C^\times$ consists of the point $z_\lambda$ in $K^+$
fixed by the action of $\lambda$ and the corresponding point
$\tau_\lambda$ in $\H^+$ is a CM point. Let
$E=\Q(\sqrt{-\n(\lambda)})$ and $\phi:E\to B$ be the embedding
determined by $\phi(\sqrt{-\n(\lambda)})=\lambda$. Then the
discriminant of this CM-point is the discriminant of the quadratic
order $R$ in $E$ such that $\phi(E)\cap\O=\phi(R)$. Note, however,
that if the CM-point $\tau_\lambda$ happens to be an elliptic point of
order $e$, then the projection $K^+\simeq\H^+\to X_0^D(N)/W_{D,N}$ is
locally $e$-to-$1$ at $\tau_\lambda$. Thus, the order of the modular
form $\psi_F(\tau)$ at $\tau_\lambda$ is $1/e$ of that of $\Psi_F(z)$ at
$z_\lambda$.

In practice, to have a simpler description of the divisor of
$\psi_F(\tau)$, we often assume that the weakly holomorphic
vector-valued modular form $F$ has the property that the only $\eta\in
L^\vee/L$ such that $c_\eta(m)\ne 0$ for some $m<0$ is $0$. In such a
case, if we assume that $\lambda$ is primitive, that is,
$\lambda/n\notin\O$ for any positive integer $n\ge 2$, then the
discriminant of the CM-point $\tau_\lambda$ is either $-\n(\lambda)$
or $-4\n(\lambda)$, depending on whether $(1+\lambda)/2$ is in $\O$ or
not. In summary, the divisor of $\psi_F(\tau)$ can be described as
follows.

\begin{Lemma} \label{lemma: divisor}
  Let $F(\tau)$, $\Psi_F(z)$, and $\psi_F(\tau)$ be as in
  the previous lemma. Assume in addition that the only $\eta\in
  L^\vee/L$ such that $c_\eta(m)\ne 0$ for some $m<0$ is $0$. Then we have
  $$
    \div\psi_F=\sum_{m<0}c_0(m)\sum_{r\in\Z^+,4m/r^2
    \text{ is a discriminant}}\frac1{e_{4m/r^2}}
    \sum_{\tau\in\CM(4m/r^2)}\tau,
  $$
  where for a negative discriminant $d$, $\CM(d)$ denotes the set of
  CM-points of discriminant $d$ (which might be empty) and $e_d$ is
  the cardinality of the stabilizer of $\tau\in\CM(d)$ in
  $N_B^+(\O)/\Q^\times$.
\end{Lemma}

We next determine when the character of a Borcherds form
$\psi_F(\tau)$ is trivial, under the assumption that the genus of
$N_B^+(\O)\backslash\H^+$ is zero.

\begin{Lemma} Assume that the genus of $X=N_B^+(\O)\backslash\H^+$ is
  zero. Let $\tau_1,\ldots,\tau_r$ be the elliptic points of $X$ and
  assume that their orders are $b_1,\ldots,b_r$, respectively. Assume
  further that, as CM-points, the discriminants of
  $\tau_1,\ldots,\tau_r$ are $d_1,\ldots,d_r$, respectively. Let
  $F(\tau)=\sum_\eta(\sum_m c_\eta(m)q^m)e_\eta$
  be a weakly holomorphic vector-valued modular form of weight $1/2$
  and type $\rho_L$ such that $O_{L,F}^+=O_L^+$ and $c_\eta(m)\in\Z$
  whenever $m\le 0$. Assume that $c_0(0)$ is even. Then the character
  associated to the modular form $\psi_F(\tau)$ is trivial
  if and only if for all $j$ such that $b_j\ne 3$, the order of
  $\Psi_F(z)$ at $z(\tau_j)$ has the same parity as $c_0(0)/2$.
\end{Lemma}

\begin{proof}
  Let $\gamma_1,\ldots,\gamma_r$ be generators of the stabilizer
  subgroups of $\tau_1,\ldots,\tau_r$ in the group
  $\Gamma=N_B^+(\O)/\Q^\times$. Since $X$ is assumed to be of
  genus zero, the group $\Gamma$ is generated by
  $\gamma_1,\ldots,\gamma_r$ with a single relation
  \begin{equation} \label{equation: Fuchsian}
    \gamma_1\ldots\gamma_r=1,
  \end{equation}
  after a suitable reindexing. (See \cite[Chapter 4]{Katok}.)

  Recall that the order of an elliptic point can only be $2$,
  $3$, $4$, or $6$. Also, an elliptic point of order $3$ or $6$ is
  necessarily a CM-point of discriminant $-3$ and a CM-point of
  discriminant $-3$ is an elliptic point of order $3$ or $6$ depending
  on whether $3\nmid DN$ or $3|DN$. In particular, an elliptic point
  of order $3$ and an elliptic point of order $6$ cannot exist at the
  same time. Moreover, on $X_0^D(N)/W_{D,N}$, there can be at most one
  CM-point of discriminant $-3$. Likewise, an elliptic point of order
  $4$ is necessarily a CM-point of discriminant $-4$ and on
  $X_0^D(N)/W_{D,N}$ there can be at most one such point. Therefore,
  there are at most two elliptic points whose orders are different
  from $2$.

  Consider the case where there is one or zero elliptic point whose
  order is different from $2$ first. By \eqref{equation: Fuchsian}, to
  show that the character $\chi$ associated to the modular form
  $\psi_F(\tau)$ is trivial, it suffices to prove that
  $\chi(\gamma_j)=1$ for $j$ with $b_j=2$.

  Observe that for $j$ with $b_j=2$, $\gamma_j$ is an element of
  order $2$ in $\Gamma$ and hence of trace zero and positive norm. Now
  by the compatibility relation \eqref{equation: compatibility}, if we
  write $\iota(\gamma_j)=\SM{c_1}{c_2}{c_3}{c_4}$ and set $k=c_0(0)$,
  then
  \begin{equation*}
    \psi_F\left(\frac{c_1\tau+c_2}{c_3\tau+c_4}\right)
   =\Psi_F\left(\frac{\n(\gamma_j)}{(c_3\tau+c_4)^2}\gamma_jz(\tau)
    \gamma_j^{-1}\right)
   =\frac{(c_3\tau+c_4)^k}{\n(\gamma_j)^{k/2}}
    \Psi_F\left(\gamma_j z(\tau)\gamma_j^{-1}\right).
  \end{equation*}
  Let $\sigma_j$ be the element of $O_L^+$ that corresponds to the
  reflection with respect to $\gamma_j$. We have
  $\sigma_j:z\to-\gamma_j z\gamma_j^{-1}$. Being a reflection,
  $\sigma_j$ acts on $\Psi_F(z)$ as $+1$ or $-1$, depending on whether
  $\Psi_F(z)$ has an even order or an odd order at the fixed point
  $z(\tau_j)$ of $\sigma_j$. Thus, assuming the order of $\Psi_F(z)$
  at $z(\tau_j)$ has the same parity as $k/2=c_0(0)/2$, we have
  \begin{equation*}
  \begin{split}
    \psi_F\left(\frac{c_1\tau+c_2}{c_3\tau+c_4}\right)
  &=\frac{(c_3\tau+c_4)^k}{\n(\gamma_j)^{k/2}}
    \Psi_F(-\sigma_jz(\tau)) \\
  &=(-1)^{k/2}\frac{(c_3\tau+c_4)^k}{\n(\gamma_j)^{k/2}}
    \Psi_F(\sigma_jz(\tau)) \\
  &=\frac{(c_3\tau+c_4)^k}{\n(\gamma_j)^{k/2}}\Psi_F(z(\tau))
   =\frac{(c_3\tau+c_4)^k}{\n(\gamma_j)^{k/2}}\psi_F(\tau).
  \end{split}
  \end{equation*}
  Therefore, if the order of $\Psi_F(z)$ at $z(\tau_j)$ has the same
  parity as $k/2=c_0(0)/2$ for all $j$ with $b_j=2$, then
  $\psi_F(\tau)$ is a modular form with trivial character on $X$.

  Now consider the remaining case where there are two elliptic points
  of order different from $2$. By the remark made earlier, the orders
  of these two elliptic points can only be $3$ and $4$ or $4$ and $6$.
  By the same argument in the previous paragraph, we find that, under
  the assumption of the lemma, for all $j$ with $b_j$ even, we have
  $\chi(\gamma_j^{b_j/2})=1$. It follows that if $b_j=4$, then
  $\chi(\gamma_j)^2=1$ and if $b_j=3$ or $b_j=6$, then
  $\chi(\gamma_j)^3=1$. Since
  $\chi(\gamma_1)\ldots\chi(\gamma_r)=1$, we conclude that
  $\chi(\gamma_j)=1$ for all $j$. This proves the lemma.
\end{proof}

For the case of $X_0^6(1)/W_6$ under consideration, there are three
elliptic points of order $2$, $4$, and $6$, respectively. They are
CM-points of discriminants $-24$, $-4$, and $-3$, respectively.
The proof of the above lemma gives us the following criterion for a
Borcherds form $\psi_F(\tau)$ to be a modular form with trivial
character on $X_0^6(1)/W_6$.

\begin{Corollary}
  \label{corollary: trivial character}
  Let $\O$ be a maximal order in the quaternion
  algebra of discriminant $6$ over $\Q$ and $L$ be the lattice formed
  by the elements of trace zero in $\O$. Suppose that
  $F(\tau)=\sum_\eta(\sum_m c_\eta(m)q^m)e_\eta$ is a weakly
  holomorphic vector-valued modular form of weight $1/2$ and type
  $\rho_L$ such that $c_\eta(m)\in\Z$ whenever $m\le 0$ and
  $O_{L,F}^+=O_L^+$. Assume in addition that
  \begin{enumerate}
  \item the only $\eta\in L^\vee/L$ such that $c_\eta(m)\ne 0$ for some
  $m<0$ is $0$, and
  \item $c_0(0)$ is even and
  $$
    \sum_{m=-r^2}c_0(m)\equiv\sum_{m=-3r^2}c_0(m)\equiv c_0(0)/2\mod 2.
  $$
  \end{enumerate}
  Then the Borcherds form $\psi_F(\tau)=\Psi_F(z(\tau))$ is a modular
  form of weight $c_0(0)$ and trivial character on the Shimura curve
  $X_0^6(1)/W_6$.
\end{Corollary}

Finally, we introduce Errthum's method for constructing $F(\tau)$
satisfying the conditions in the lemma above \cite{Errthum}.
Here we consider general Eichler orders in an
indefinite quaternion algebra over $\Q$.

The first lemma shows that we can construct $F(\tau)$ out of a
scalar-valued modular form with suitable properties. To state the
required properties, we let $\chi_\theta$ denote the character
associated to the Jacobi theta function
$\theta(\tau)=\sum_{n\in\Z}q^{n^2}$.
That is, $\chi_\theta$ is the character satisfying
$$
  \theta(\gamma\tau)=\chi_\theta(\gamma)(c\tau+d)^{1/2}\theta(\tau)
$$
for all $\gamma=\SM abcd\in\Gamma_0(4)$ and all $\tau\in\H^+$. For a
scalar-valued modular form $f(\tau)$ of weight $k\in\frac12\Z$ on
$\Gamma_0(M)$ and $\gamma=\SM abcd\in\Gamma_0(M)$, we let
$$
  f\big|_\gamma(\tau)=(c\tau+d)^{-k}f(\gamma\tau).
$$
We observe that the level $M$ of the lattice $L$ is always a multiple
of $4$ for any $D$ and $N$.

\begin{Lemma}[{\cite[Theorem 4.2.9]{Barnard}}]
  \label{lemma: Barnard}
  Let $M$ be the level of the lattice $L$. Suppose that $f(\tau)$ is a
  weakly holomorphic scalar-valued modular form of weight $1/2$ such
  that
  $$
    f(\gamma\tau)=\chi_\theta(\gamma)(c\tau+d)^{1/2}f(\tau)
  $$
  for all $\gamma=\SM abcd\in\Gamma_0(M)$. Then the function
  $F_f(\tau)$ defined by
  \begin{equation} \label{equation: Ff}
    F_f(\tau)=\sum_{\gamma\in\wt{\Gamma}_0(M)\backslash
    \wt\SL(2,\Z)}f\big|_\gamma(\tau)\rho_L(\gamma^{-1})e_0
  \end{equation}
  is a weakly holomorphic vector-valued modular form of weight $1/2$
  and type $\rho_L$.
\end{Lemma}

\begin{Lemma}[{\cite[Proposition 5.4]{Errthum}}]
  \label{lemma: infinity cusp}
  Suppose that the weakly holomorphic modular form $f(\tau)$ in the
  above lemma has a pole only at the infinity cusp. Then the Fourier
  expansion $F_f(\tau)=\sum_\eta(\sum_m c_\eta(m)q^m)e_\eta$ satisfies
  $c_\eta(m)=0$ whenever $\eta\ne 0$ and $m<0$.
\end{Lemma}

\begin{Lemma}[{\cite[Theorem 5.8]{Errthum}}] Let $f(\tau)$ and
  $F_f(\tau)$ be given as in the previous lemmas. Then for any
  $\eta,\eta'\in L^\vee$ with $\gen{\eta,\eta}=\gen{\eta',\eta'}$, the
  $e_\eta$-component and the $e_{\eta'}$-component of $F_f(\tau)$ are
  equal. Consequently, we have $O_{L,F_f}^+=O_L^+$.
\end{Lemma}

It remains to construct scalar-valued modular forms $f(\tau)$
satisfying the condition in Lemma \ref{lemma: Barnard}.

\begin{Lemma}[{\cite[Theorem 6.2]{Borcherds-Duke}}]
  \label{lemma: eta}
  Let $M$ be the level of the lattice $L$. Suppose that $r_d$, $d|M$,
  are integers satisfying the conditions
  \begin{enumerate}
  \item $\sum_{d|M} r_d=1$,
  \item $|L^\vee/L|\prod_{d|M} d^{r_d}$ is a square in $\Q^\times$,
  \item $\sum_{d|M} dr_d\equiv 0\mod 24$, and
  \item $\sum_{d|M} (M/d)r_d\equiv 0\mod 24$.
  \end{enumerate}
  Then $\prod_{d|M}\eta(d\tau)^{r_d}$ is a weakly holomorphic modular form
  satisfying the condition for $f(\tau)$ in Lemma \ref{lemma: Barnard}.
\end{Lemma}

We now consider the case of $X_0^6(1)/W_6$.

\begin{Proposition} \label{proposition: realization 2}
Consider the case $X_0^6(1)/W_6$. Let
$$
  f(\tau)=2\frac{\eta(2\tau)\eta(3\tau)^2\eta(4\tau)^4\eta(6\tau)^4}
  {\eta(12\tau)^{10}}+2\frac{\eta(\tau)\eta(2\tau)^3\eta(6\tau)^2}
  {\eta(3\tau)\eta(4\tau)\eta(12\tau)^3}.
$$
and
$$
  g(\tau)=2\frac{\eta(\tau)\eta(2\tau)^3\eta(6\tau)^2}
  {\eta(3\tau)\eta(4\tau)\eta(12\tau)^3}
$$
Let $F_f(\tau)$ and $F_g(\tau)$ be defined as in \eqref{equation: Ff}.
Then $\psi_{F_f}(\tau)$ and $\psi_{F_g}(\tau)$ span the
one-dimensional spaces of holomorphic modular form on $X_0^6(1)/W_6$ of
weight $8$ and $12$, respectively.
\end{Proposition}

\begin{proof} The two eta-products were found in \cite[Page
  848]{Errthum}. Here we give a quick explanation.

  In the case of $X_0^6(1)/W_6$, the lattice $L$ has level
  $12$ and $|L^\vee/L|=72$. The two eta-products clearly
  satisfy the four conditions in Lemma \ref{lemma: eta}. Now the
  congruence subgroup $\Gamma_0(12)$ has $6$ cusps, represented by
  $1/c$, $c|12$. The orders of the eta functions $\eta(d\tau)$ at
  these cusps, multiplied by $24$, are given by the table.
  $$ \extrarowheight3pt
  \begin{array}{c|cccccc} \hline\hline
    & 1/1 & 1/2 & 1/4 & 1/3 & 1/6 & 1/12 \\ \hline
  \eta(\tau) & 12 & 3 & 3 & 4 & 1 & 1 \\
  \eta(2\tau)&  6 & 6 & 6 & 2 & 2 & 2 \\
  \eta(4\tau)&  3 & 3 &12 & 1 & 1 & 4 \\
  \eta(3\tau)&  4 & 1 & 1 &12 & 3 & 3 \\
  \eta(6\tau)&  2 & 2 & 2 & 6 & 6 & 6 \\
  \eta(12\tau)& 1 & 1 & 4 & 3 & 3 &12 \\ 
  \hline\hline
  \end{array}
  $$
  From the table, we see that the two eta-products have only a pole at
  the cusp $1/12\sim\infty$. Thus, by Lemma \ref{lemma: infinity
    cusp}, the divisors of $\psi_{F_f}(\tau)$ and $\psi_{F_g}(\tau)$
  are determined by the $e_0$-components of the Fourier expansions of
  $F_f(\tau)$ and $F_g(\tau)$. Since
  $f(\tau)=2q^{-3}-6-18q+\cdots$ and $g(\tau)=2q^{-1}-2-8q+8q^2+\cdots$,
  the $e_0$-components of $F_f(\tau)$ and $F_g(\tau)$ are
  $$
    2q^{-3}+c_0+\cdots, \qquad 2q^{-1}+d_0+\cdots
  $$
  for some $c_0$ and $d_0$, respectively. The numbers $c_0$ and $d_0$
  are complicated to compute directly from the definition of $F_f$ and
  $F_g$. Here we observe that, by Lemma \ref{lemma: divisor},
  $$
    \div\psi_{F_f}(\tau)=\frac13P_{-3}, \qquad
    \div\psi_{F_g}(\tau)=\frac12P_{-4},
  $$
  where $P_{-3}$ and $P_{-4}$ denote the unique CM-points of discriminants
  $-3$ and $-4$, respectively. (Note that there does not exist a
  CM-point of discriminant $-12$ on $X_0^6(1)/W_6$.) Therefore, the
  weight of $\psi_{F_f}(\tau)$ must be $8$ and the weight of
  $\psi_{F_g}(\tau)$ must be $12$. In other words, we have $c_0=8$ and
  $d_0=12$. Then, by Corollary \ref{corollary: trivial character},
  $\psi_{F_f}(\tau)$ and $\psi_{G_f}(\tau)$ must be modular forms on
  $X_0^6(1)/W_6$ with trivial characters. This proves the proposition.
\end{proof}

Combining Proposition \ref{proposition: realization 1} and Proposition
\ref{proposition: realization 2}, we find that
\begin{equation} \label{equation: C1}
  \psi_{F_f}(\tau)=C_1\left({}_2F_1\left(\frac1{24},\frac5{24};\frac34;s
  \right)+\frac1{\sqrt[4]{12}\omega_{-4}^2}s^{1/4}
  {}_2F_1\left(\frac7{24},\frac{11}{24};\frac54;s\right)\right)^8
\end{equation}
and
\begin{equation} \label{equation: C2}
  \psi_{F_g}(\tau)=C_2\left({}_2F_1\left(\frac1{24},\frac7{24};\frac56;t
  \right)-\frac{e^{-2\pi i/8}}{\sqrt[6]2\omega_{-3}^2}t^{1/6}
  {}_2F_1\left(\frac5{24},\frac{11}{24};\frac76;t\right)\right)^{12}
\end{equation}
for some complex numbers $C_1$ and $C_2$. To determine the absolute
values of these two numbers, we shall use Schofer's formula for values
of Borcherds forms at CM-points.
\end{section}

\begin{section}{Schofer's formula for values of Borcherds forms at
    CM-points}
\label{section: Schofer}
Let $\O$ be an Eichler order of level $N$ in an indefinite quaternion
algebra of discriminant $D$ over $\Q$. Throughout this section, we
assume that the level $N$ is squarefree and the symbol $d$ always
denote a negative fundamental discriminant. Let
$L=\{\alpha\in\O:\tr(\alpha)=0\}$ be the lattice of signature $(1,2)$
formed by the elements of trace $0$ in $\O$. We retain all the
notations $\gen{\cdot,\cdot}$, $V(\Q)$, $V(\R)$, $V(\C)$, $K$,
$\widetilde K^+$, $O_L$, $O_{L,F}$, and etc. used in the previous
section. Here let us summarize Schofer's formula \cite{Schofer} for
average values of Borcherds forms at CM-points first. The explanation
of the terms involved will be given later.

\begin{theorem}[{\cite[Corollaries 1.2 and 3.5]{Schofer}}]
  \label{theorem: Schofer}
  Let $F(\tau)=\sum_\eta(\sum_m c_\eta(m)q^m)e_\eta$ be a weakly
  holomorphic vector-valued modular form of weight $1/2$ and type
  $\rho_L$ such that $O_{L,F}^+=O_L^+$ and $c_\eta(m)\in\Z$ whenever
  $m\le 0$. Let $\Psi_F(z)$ be the Borcherds form associated $F(\tau)$
  and $\psi_F(\tau)=\Psi_F(z(\tau))$ be the modular form of weight
  $c_0(0)$ on $X_0^D(N)/W_{D,N}$ as described in Lemma \ref{lemma:
    psi}, where $z(\tau)$ is given by \eqref{equation: z(tau)}. Let
  $d<0$ be a fundamental discriminant such that the set $\CM(d)$ of
  CM-points of discriminant $d$ is not empty and that the support of
  $\div\psi(\tau)$ does not intersect $\CM(d)$. Then we have
  \begin{equation*}
  \begin{split}
   &\sum_{\tau\in\CM(d)}\log\left|\psi_F(\tau)(\Im\tau)^{c_0(0)/2}\right| \\
   &\qquad\qquad=-\frac{|\CM(d)|}4\left(\sum_{\eta\in L^\vee/L}
    \sum_{m\ge 0}c_\eta(-m)\kappa_\eta(m)
    +c_0(0)(\Gamma'(1)+\log(2\pi))\right).
  \end{split}
  \end{equation*}
\end{theorem}

\begin{Remark}
\begin{enumerate}
\item Note that the formula given in \cite{Schofer} is valid for
  Borcherds forms associated to lattices of general signature $(n,2)$.
  Here we have specialized the formulas to the cases under our
  consideration. Note also that in \cite{Schofer}, the left-hand side
  of the formula has $\Psi_F(z)|y|^{c_0(0)/2}$ in place of
  $\psi_F(\tau)(\Im\tau)^{c_0(0)/2}$, where $z=x+iy\in\widetilde K^+$
  and $|y|=\sqrt{|\gen{y,y}|}$. By a direct computation, we find that
  for $z=z(\tau)$ given in \eqref{equation: z(tau)}, we have
  $|y|=\Im\tau$. Notice that in general, for any modular form
  $\psi(\tau)$ of weight $k$ on a Fuchsian subgroup $\Gamma$ of
  $\SL(2,\R)$, we have
  $|\psi(\gamma\tau)(\Im\gamma\tau)^{k/2}|=|\psi(\tau)(\Im\tau)^{k/2}|$
  for any $\tau\in\H^+$ and $\gamma\in\Gamma$. Thus, the left-hand
  side of the formula does not depend on the choice of representatives
  of the CM-points.
\item Let $\chi_d$ be the Kronecker character associated to the field
  $\Q(\sqrt d)$ and
  $$
    \Lambda(s,\chi_d)=\left(\frac{\pi}{|d|}\right)^{-(1+s)/2}
    \Gamma\left(\frac{1+s}2\right)L(s,\chi_d)
  $$
  be the complete $L$-function associated to $\chi_d$. In
  \cite{Schofer}, the term $\kappa_0(0)$ was given as
  $$
    \kappa_0(0)=2\frac{\Lambda'(1,\chi_d)}{\Lambda(1,\chi_d)}
  $$
  under a certain assumption. (Note that our definition of
  $\Lambda(s,\chi_d)$ is different from that in \cite{Schofer}.) Later
  on, we will prove that for the cases under our consideration, we have
  $$
    \kappa_0(0)=2\frac{\Lambda'(1,\chi_d)}{\Lambda(1,\chi_d)}
   +\sum_{p|D/(D,d)}\frac{p-1}{p+1}\log p
   +\sum_{p|N/(N,d)}\log p,
  $$
  where the last two summations run over all prime divisors $p$ of
  $D/(D,d)$ and $N/(N,d)$, respectively. 
\item From the functional equation
  $\Lambda(s,\chi_d)=\Lambda(1-s,\chi_d)$ for $\Lambda(s,\chi_d)$, we
  deduce that
  \begin{equation} \label{equation: Lambda'}
    2\frac{\Lambda'(1,\chi_d)}{\Lambda(1,\chi_d)}
   =\log\frac{4\pi}{|d|}-\Gamma'(1)-2\frac{L'(0,\chi_d)}{L(0,\chi_d)}.
  \end{equation}
  By the Chowla-Selberg formula, we have
  $$
    e^{L'(0,\chi_d)/2L(0,\chi_d)}=\frac1{\sqrt{|d|}}\prod_{a=1}^{|d|-1}
    \Gamma\left(\frac a{|d|}\right)^{\chi_d(a)\mu_d/4h_d}=\omega_d.
  $$
  This shows that the value of a suitable modular form of
  weight $k$ on $X_0^D(N)/W_{D,N}$ at a CM-point of discriminant $d$
  will be an algebraic multiple of $\omega_d^k$, agreeing with the results
  of \cite{Shimura-periods} and \cite{Yoshida}.
\end{enumerate}
\end{Remark}

We now explain what the terms $\kappa_\eta(m)$ are.
Recall that each CM-point $\tau$ of discriminant $d$ corresponds to an
embedding $\phi:\Q(\sqrt d)\to B$ such that $\phi(\Q(\sqrt
d))\cap\O=\phi(R_d)$, where $R_d$ is the imaginary quadratic order of
discriminant $d$. To be more precise, $\tau$ is the common fixed point
of $\phi(R_d)$ in the upper half-plane. Let $\lambda=\phi(\sqrt d)$.
Then $\lambda$ is an element of positive norm in $L$ and the set
$U=\lambda^\perp=\{\alpha\in V(\Q):\gen{\lambda,\alpha}=0\}$ is a
negative $2$-plane isomorphic to $\Q(\sqrt d)$ in the sense that there
is an isomorphism $h:U\to\Q(\sqrt d)$ as vector spaces over $\Q$ and a
negative rational number $c$ such that
$c\gen{\alpha,\beta}=\tr_\Q^{\Q(\sqrt
  d)}(h(\alpha)\overline{h(\beta)})$ for all $\alpha,\beta\in U$.

Let $L_+=L\cap\Q\lambda$ and $L_-=L\cap U$. We have
$$
  L_++L_-\subset L\subset L^\vee\subset L_+^\vee+L_-^\vee.
$$
For $\mu\in L_-^\vee/L_-$, let $\varphi_\mu:U\to\C$ be the
characteristic function of $\mu+L_-$. Then for each $\mu\in
L_-^\vee/L_-$, we have the incoherent Eisenstein series
$E(\tau,s;\varphi_\mu,1)$ of weight $1$
\cite{Kudla-Rapoport-Yang,Kudla-Rapoport-Yang2,Kudla-Yang,
Schofer}. Write
$\tau=u+iv$ and let
$$
  E(\tau,s;\varphi_\mu,1)=\sum_m A_\mu(s,m,v)q^m, \quad q=e^{2\pi i\tau},
$$
be the Fourier expansion of $E(\tau,s;\varphi_\mu,1)$. The Eisenstein
series $E(\tau,s;\varphi_\mu,1)$ vanishes at $s=0$. Thus,
$$
  A_\mu(s,m,v)=b_\mu(m,v)s+O(s^2)
$$
near $s=0$ for some function $b_\mu(m,v)$. We define
\begin{equation} \label{equation: kappa-}
  \kappa_\mu^-(m):=\begin{cases}
  \displaystyle\lim_{v\to\infty} b_\mu(m,v), &\text{if }m>0, \\
  \displaystyle\lim_{v\to\infty}(b_0(0,v)-\log v), &\text{if }m=0 
    \text{ and }\mu=0, \\
  0, &\text{else}. \end{cases}
\end{equation}
Then the term $\kappa_\eta(m)$ in Schofer's formula is defined by
\begin{equation} \label{equation: kappa}
  \kappa_\eta(m)=\sum_{\mu\in L/(L_++L_-)}\sum_{x\in\eta_++\mu_++L_+}
  \kappa_{\eta_-+\mu_-}^-(m-\gen{x,x}/2),
\end{equation}
where for $\mu\in L/(L_++L_-)$ and $\eta\in L^\vee/L$, we write
$\mu=\mu_++\mu_-$ and $\eta=\eta_++\eta_-$ with
$\mu_+,\eta_+\in\Q\lambda$ and $\mu_-,\eta_-\in U$. The terms
$\kappa_\eta(m)$ look very complicated, but nonetheless are computable
using the fact that $A_\mu(s,m,v)q^m$ can be written as a product of
local Whittaker functions
\cite{Kudla-Rapoport-Yang,Kudla-Rapoport-Yang2,Kudla-Yang}, which can
be computed using formulas in \cite{Kudla-Yang,Yang-density}. Here we
briefly describe a general strategy to compute $A_\mu(s,m,v)$ and
$\kappa_\mu^-(m)$ efficiently, following \cite{Errthum,Kudla-Yang}.

In general, for $\mu\in L_-^\vee/L_-$, we have $A_\mu(s,m,v)=0$ unless
$\gen{\mu,\mu}/2-m\in\Z$ and when $\gen{\mu,\mu}/2-m\in\Z$ holds, we
have
$$
  A_\mu(s,m,v)q^m=\delta_{\mu,m}v^{s/2}+W_{m,\infty}(\tau,s)\prod_{p<\infty}
  W_{m,p}(s,\varphi_{\mu,p}),
$$
where
$$
  \delta_{\mu,m}=\begin{cases}
  1, &\text{if }\mu=0\text{ and }m=0, \\
  0, &\text{else}, \end{cases}
$$
and $W_{m,\infty}(\tau,s)$ and $W_{m,p}(s,\varphi_{\mu,p})$ are
the local Whittaker factors at $\infty$ and $p$, respectively.
(See \cite[Section 2]{Kudla-Rapoport-Yang}.) Let $\Delta$ be the
discriminant of the lattice $L_-$. When a prime $p$ does not divide
$\Delta$ and the $p$-adic valuation $v_p(m)$ is zero, Equation
(4.4) and Theorems 4.3 and 4.4 of \cite{Kudla-Yang} yield
$$
  W_{m,p}(s,\varphi_{\mu,p})=\gamma_p(1-\chi_d(p)p^{-1-s}),
$$
where $\gamma_\infty$ and $\gamma_p$ are certain explicit constants
that do not have any effect on the calculation since
$\gamma_\infty\prod_p\gamma_p=1$. Therefore, assuming $m>0$, letting
\begin{equation} \label{equation: S}
  S_{m,\mu}=\{p:~p|\Delta \text{ or }v_p(m)>0\},
\end{equation}
and using the formula for $W_{m,\infty}$ in Proposition 2.3 of
\cite{Kudla-Yang}, we find
$$
  A_\mu(m,s,v)=-\frac{2\pi}{L(s+1,\chi_d)}\prod_{p\in S_{m,\mu}}
  \frac{W_{m,p}(s,\varphi_{\mu,p})}{1-\chi_d(p)p^{-1-s}}.
$$
As $A_\mu(m,0,v)=0$, there exists at least a prime $p'$ in $S_{m,\mu}$
such that $W_{m,p'}(0,\varphi_{\mu,p})$. Taking the derivative of the
above expression and evaluating at $s=0$, we obtain the following lemma.

\begin{Lemma} \label{lemma: kappa(m)}
  Assume that $m>0$ and let all the notations be given as
  in the discussion. We have
  $$
  \kappa_\mu^-(m)=-\frac{\mu_d\sqrt{|d|}}{h_d}
  \frac{W'_{m,p'}(0,\varphi_{\mu,p'})}{1-\chi_d(p')/p'}
  \prod_{p\in S_{m,\mu},p\ne p'}\frac{W_{m,p}(0,\varphi_{\mu,p})}{1-\chi_d(p)/p},
  $$
  where $\mu_d$ and $h_d$ denote the number of roots of unity and the
  class number of $\Q(\sqrt d)$, respectively.
\end{Lemma}

This is essentially Theorem 6.3 of \cite{Errthum}. We now consider the
term $\kappa_0^-(0)$. If the discriminant $\Delta$ of $L_-$ is precisely
$|d|$, then, again, Theorems 4.3 and 4.4 of \cite{Kudla-Yang} show
that
$$
  W_{0,p}(s,\varphi_{0,p})=\gamma_p\frac{1-\chi_d(p)p^{-1-s}}{1-\chi_d(p)p^{-s}}
$$
so that
$$
  A_0(s,0,v)=v^{s/2}-v^{-s/2}\frac{\Lambda(s,\chi_d)}{\Lambda(s+1,\chi_d)}
$$
and
$$
  b_0(0,v)=\log v+\frac{\Lambda'(1,\chi_d)}{\Lambda(1,\chi_d)}
  -\frac{\Lambda'(0,\chi_d)}{\Lambda(0,\chi_d)}
  =\log v+2\frac{\Lambda'(1,\chi_d)}{\Lambda(1,\chi_d)}.
$$
(See \cite[Lemma 2.20]{Schofer}.)

In general, the discriminant $\Delta$ of $L_-$ may not be exactly
$|d|$. Let
$$
  S=\{p:~p|(\Delta/d)\}.
$$
Then we have
\begin{equation} \label{equation: A0}
  A_0(s,0,v)=v^{s/2}-v^{-s/2}\frac{\Lambda(s,\chi_d)}{\Lambda(s+1,\chi_d)}
  \prod_{p\in S}\frac{(1-\chi_d(p)p^{-s})W_{0,p}(s,\varphi_{0,p})}{1-\chi_d(p)p^{-1-s}}.
\end{equation}
Let $G(s)$ denote the product on the right. Since $A_0(0,0,v)$ is
identically $0$, we must have $G(0)=1$. From this, we deduce the
following lemma.

\begin{Lemma} \label{lemma: kappa(0)}
  Let all the notations be given as above. We have
  $$
    \kappa_0^-(0)=2\frac{\Lambda'(1,\chi_d)}{\Lambda(1,\chi_d)}
   -\frac d{ds}\log G(s)\Big|_{s=0}.
  $$
\end{Lemma}

We now determine $G(s)$ for the cases under our consideration.
In the following lemma, for a prime $p$, we let $L_p=L_-\otimes_\Z\Z_p$.

\begin{Lemma} \label{lemma: local lattices}
  Let all the notations be given as in the preceding
  discussion. Assume that the level $N$ of the Eichler order $\O$ is
  squarefree and that $d<0$ is a fundamental discriminant.
  \begin{enumerate}
  \item Let $p$ be an odd prime. Then there exists a basis
    $\{\ell_1,\ell_2\}$ for $L_p$ and
    $\epsilon_1,\epsilon_2\in\Z_p$ with
    $\epsilon_1\epsilon_2=-d$ such that the Gram matrix
    $(\gen{\ell_i,\ell_j})$ is equal to
    $$
      \begin{cases}
      \M{\epsilon_1}00{\epsilon_2}, &\text{if }p|(DN,d)\text{ or if }
      p\nmid DN, \\
      p\M{\epsilon_1}00{\epsilon_2}, &\text{if }p|DN\text{ but }
      p\nmid d. \end{cases}
    $$
  \item Assume that $d\equiv 0\mod 4$. Then there exists a basis
    $\{\ell_1,\ell_2\}$ for $L_2$ and $\epsilon_1,\epsilon_2\in\Z_2$
    with $\epsilon_1\epsilon_2=-d/4$ such that the Gram matrix is
    $$
      2\M{\epsilon_1}00{\epsilon_2}.
    $$
  \item Assume that $d\equiv 1\mod 8$ (and $2\nmid D$). Then there
    exists a basis $\{\ell_1,\ell_2\}$ for $L_2$ and
    $\epsilon\in\Z_2^\times$ such that the Gram matrix is
    $$
      \begin{cases}
      2\epsilon\M0110, &\text{if }2|N, \\
      \epsilon\M0110, &\text{if }2\nmid N. \end{cases}
    $$
  \item Assume that $d\equiv 5\mod 8$ (and $2\nmid N$). Then there
    exists a basis $\{\ell_1,\ell_2\}$ for $L_2$ and
    $\epsilon\in\Z_2^\times$ such that the Gram matrix is
    $$
      \begin{cases}
      2\epsilon\M2112 &\text{if }2|D, \\
      \epsilon\M2112
      &\text{if }2\nmid D, \end{cases}
    $$
  \end{enumerate}
\end{Lemma}

\begin{proof} Assume that $p$ is an odd prime. Consider the case when
  $p$ divides $DN$ first. There exists a basis $\{e_1,e_2,e_3\}$
  for $L\otimes\Z_p$ such that
  $$
    (\gen{e_i,e_j})=\begin{pmatrix}
    2\mu_1 & 0 & 0 \\ 0 & 2\mu_2p & 0 \\ 0 & 0 & 2\mu_1\mu_2p\end{pmatrix},
  $$
  where $\mu_1$ and $\mu_2$ are elements in $\Z_p^\times$ with the
  property that the Hilbert symbol $(-\mu_1,-\mu_2p)_p$ is $1$ or $-1$
  depending on whether $p|N$ or $p|D$.

  Assume that $\lambda=c_1e_1+c_2e_2+c_3e_3$. If $p|d$, then we have
  $p|c_1$ and at least one of $c_2$ and $c_3$ must be in
  $\Z_p^\times$. Without loss of generality, we assume that
  $c_2\in\Z_p^\times$. Then $L_p$ is spanned by
  $-c_2\mu_2e_1+(c_1/p)\mu_1e_2$ and $c_3\mu_1e_2-c_2e_3$. The Gram
  matrix of $L_p$ with respect to this basis has determinant
  $-(2c_2\mu_1\mu_2)^2d$. It follows that there is a basis
  $\{\ell_1,\ell_2\}$ for $L_p$ such that
  $(\gen{\ell_i,\ell_j})=\SM{\epsilon_1}00{\epsilon_2}$ with the
  properties $\epsilon_1,\epsilon_2\in\Z_p$ and
  $\epsilon_1\epsilon_2=-d$.

  If $p\nmid d$, then $p\nmid c_1$. We find that $L_p$ is spanned by
  $-c_2\mu_2pe_2+c_1\mu_1e_2$ and $-c_3\mu_2p+c_1e_3$. The Gram matrix
  of $L_p$ with respect to this basis is inside $pM(2,\Z_p)$ and its
  determinant is $-(2c_1\mu_1\mu_2p)^2d$. It follows that there exists
  a basis $\{\ell_1,\ell_2\}$ for $L_p$ such that
  $(\gen{\ell_i,\ell_j})=\SM{\epsilon_1p}00{\epsilon_2p}$ with
  $\epsilon_1,\epsilon_2\in\Z_p$ and $\epsilon_1\epsilon_2=-d$. The
  proof of the case $p\nmid DN$ is similar and is omitted.

  Now consider the case $p=2$. If $2\nmid DN$, then $\O\otimes_\Z\Z_2$
  is isomorphic to $M(2,\Z_2)$. Thus, we may assume that $L\otimes_\Z\Z_2$
  is $\{\alpha\in M(2,\Z_2):\tr(\alpha)=0\}$ so that $e_1=\SM100{-1}$,
  $e_2=\SM0100$, and $e_3=\SM0010$ form a basis for $L\otimes_\Z\Z_2$. Let
  $c_1,c_2,c_3$ be the elements in $\Z_2$ such that
  $$
    c_1e_1+c_2e_2+c_3e_3=\begin{cases}
    \lambda, &\text{if }d\equiv 1\mod 4, \\
    \lambda/2, &\text{if }d\equiv 0\mod 4. \end{cases}
  $$
  When $d\equiv 1\mod 4$, the element $\lambda$ satisfies
  $(1+\lambda)/2\in\O\otimes\Z_2$, which implies that $2\nmid c_1$ and
  $2|c_2,c_3$. Therefore, the lattice $L_2=L_-\otimes_\Z\Z_2$ is
  spanned by $-(c_2/2)e_1+c_1e_3$ and $-(c_3/2)e_1+c_1e_2$. The Gram
  matrix relative to this basis is
  $$
    \M{-c_2^2/2}{-c_1^2-c_2c_3/2}{-c_1^2-c_2c_3/2}{-c_3^2/2}
  $$
  with determinant $-c_1^2(c_1^2+c_2c_3)\equiv-d\mod 8$.
  By Lemma 8.4.1 of \cite{Cassels}, there is a basis
  $\{\ell_1,\ell_2\}$ for $L_2$ and $\epsilon\in\Z_2^\times$ such that
  the Gram matrix is
  $$
    \begin{cases}
    \epsilon\M0110, &\text{if }d\equiv 1\mod 8, \\
    \epsilon\M2112, &\text{if }d\equiv 5\mod 8. \end{cases}
  $$

  If $d\equiv 0\mod 4$, then $c_2$ and $c_3$ cannot be both even since
  $-c_1^2-c_2c_3=-d/4\equiv 1,2\mod 4$. Assume that $2\nmid c_2$. Then
  $L_2$ is spanned by $c_2e_1-2c_1e_3$ and $c_2e_2-c_3e_3$. The Gram
  matrix with respect to this basis is
  $$
    \M{-2c_2^2}{2c_1c_2}{2c_1c_2}{2c_2c_3}.
  $$
  It follows from Lemma 8.4.1 of \cite{Cassels} that there exists a
  basis $\{\ell_1,\ell_2\}$ for $L_2$ such that the Gram matrix is
  $$
    2\M{\epsilon_1}00{\epsilon_2}
  $$
  with $\epsilon_1,\epsilon_2\in\Z_2$ and $\epsilon_1\epsilon_2=-d/4$.
  This proves the case $2\nmid DN$.

  The proof of the case $2|DN$ is similar. We remark that when $2|N$,
  we have $\O\otimes_\Z\Z_2\simeq \SM{\Z_2}{\Z_2}{2\Z_2}{\Z_2}$ and
  when $2|D$, we have $B\otimes\Q_2\simeq\JS{-1,-1}{\Q_2}$ and the
  maximal order in $\JS{-1,-1}{\Q_2}$ is
  $\Z_2+\Z_2I+\Z_2J+\Z_2(1+I+J+IJ)/2$. The rest of proof is similar to
  that in the other cases and is omitted.
\end{proof}

\begin{Corollary} \label{corollary: lattice}
  Let all the notations and assumptions be given as before. Let
  $$
    r=\prod_{p|DN/(DN,d)}p.
  $$
  Then the Gram matrix of $L_-$ is equivalent to $-rM$ for some
  positive definite integral matrix $M$ of determinant $|d|$. In
  particular, the discriminant of $L_-$ is $r^2|d|$.
\end{Corollary}

\begin{Lemma} \label{lemma: kappa0(0)}
  Assume that $N$ is squarefree and $d<0$ is a
  fundamental discriminant. Let $\chi_d$, $\Lambda(s,\chi_d)$,
  $\lambda$, $L^+$, and $L_-$ be defined as above. Let
  $\kappa_\mu^-(m)$ and $\kappa_\eta(m)$ be defined as in
  \eqref{equation: kappa-} and \eqref{equation: kappa},
  respectively. Then we have
  $$
    \kappa_0^-(0)=\kappa_0(0)=2\frac{\Lambda'(1,\chi_d)}{\Lambda(1,\chi_d)}
   +\sum_{p|D/(D,d)}\frac{p-1}{p+1}\log p+\sum_{p|N/(N,d)}\log p,
  $$
  where the two sums run over prime divisors of $D/(D,d)$ and
  $N/(N,d)$, respectively.
\end{Lemma}

\begin{proof}
  Consider the case when an odd prime $p$ divides $DN/(DN,d)$, i.e.,
  $p|DN$, but $p\nmid d$. By Lemma \ref{lemma: local lattices},
  the Gram matrix of $L_p=L_-\otimes_\Z\Z_p$ is equivalent to
  $p\SM{\epsilon_1}00{\epsilon_2}$ for some
  $\epsilon_1,\epsilon_2\in\Z_p$ with $\epsilon_1\epsilon_2=-d$. We
  shall apply Theorem 4.3 of \cite{Kudla-Yang} with $\mu=0$ and $m=0$.
  Using the notations in Section 4.2 of \cite{Kudla-Yang}, we have
  $H_\mu=\{1,2\}$, $K_0(\mu)=\infty$,
  $$
    L_\mu(k)=\begin{cases}
    \{1,2\}, &\text{if }k\text{ is even}, \\
    \emptyset, &\text{if }k\text{ is odd}, \end{cases}
  $$
  $d_\mu(k)=1$ for all $k$, $\epsilon_\mu(k)=\chi_d(p)^{k-1}$,
  $t_\mu(m)=0$, and $a_\mu(m)=\infty$. Thus, the combination of (4.4)
  and Theorem 4.3 of \cite{Kudla-Yang} yields
  $$
    \frac{W_{0,p}(s,\varphi_{0,p})}{\gamma_pp^{-1}}
   =1+\left(1-\frac1p\right)\sum_{k=1}^\infty p\chi_d(p)^{k-1}p^{-ks}
   =1+(p-1)\frac{p^{-s}}{1-\chi_d(p)p^{-s}}.
  $$
  That is
  \begin{equation} \label{equation: Whittaker p}
    W_{0,p}(s,\varphi_{0,p})
  =\gamma_p\cdot\frac{1-\chi_d(p)p^{-1-s}}{1-\chi_d(p)p^{-s}}\cdot
    \frac{1+(p-1-\chi_d(p))p^{-s}}{p-\chi_d(p)p^{-s}}.
  \end{equation}

  For the case $2|DN/(DN,d)$, i.e., $2|DN$ and $d\equiv 1\mod 4$, we
  use the results in Section 4.3 of \cite{Kudla-Yang}. Consider the
  case $d\equiv 1\mod 8$ and $2|N$ first. By Lemma
  \ref{lemma: local lattices}, the Gram matrix of $L_2$ is equivalent
  to $2\epsilon\SM0110$. Following the notations in Section 4.3 of
  \cite{Kudla-Yang}, we have $H_\mu=N_\mu=\emptyset$, $M_\mu=\{1\}$,
  $L_\mu(k)=\emptyset$,
  $d_\mu(k)=p_\mu(k)=\epsilon_\mu(k)=\delta_\mu(k)=1$ for all $k\ge 1$,
  $K_0(\mu)=\infty$, and $t_\mu=\nu=0$. Thus, Theorem 4.4 and (4.4)
  of \cite{Kudla-Yang} yield
  \begin{equation} \label{equation: Whittaker 2 1}
    W_{0,2}(s,\varphi_{0,p})=\frac{\gamma_2}2\left(1+
    2^{-s}+2^{-2s}+\cdots\right)
   =\gamma_2\cdot\frac{1-2^{-1-s}}{1-2^{-s}}\cdot
    \frac{1}{2-2^{-s}}.
  \end{equation}
  For the case $d\equiv 5\mod 8$ and $2|D$, Lemma \ref{lemma: local
    lattices} shows that the Gram matrix is equivalent to
  $2\epsilon\SM2112$ for some $\epsilon\in\Z_2^\times$. In this case,
  we have $H_\mu=M_\mu=\emptyset$, $N_\mu=\{1\}$,
  $L_\mu(k)=\emptyset$, $d_\mu(k)=\epsilon_\mu(k)=\delta_\mu(k)=1$ for
  $k\ge 1$, $p_\mu(k)=(-1)^{k-1}$, $K_0(\mu)=\infty$, and
  $t_\mu=\nu=0$. Then Theorem 4.4 of \cite{Kudla-Yang} shows that
  \begin{equation} \label{equation: Whittaker 2 2}
    W_{0,2}(s,\varphi_{0,p})=\frac{\gamma_2}2\left(
    1+2^{-s}-2^{-2s}+2^{-3s}-\cdots\right)
  =\gamma_2\cdot\frac{1+2^{-1-s}}{1+2^{-s}}\cdot\frac{1+2^{1-s}}{2+2^{-s}}.
  \end{equation}
  From \eqref{equation: A0}, \eqref{equation: Whittaker p},
  \eqref{equation: Whittaker 2 1}, and \eqref{equation: Whittaker 2
    2}, we see that
  \begin{equation*}
  \begin{split}
    A_0(s,0,v)
  &=v^{s/2}-v^{-s/2}\frac{\Lambda(s,\chi_d)}
      {\Lambda(s+1,\chi_d)}\prod_{p|D/(D,d)}\frac{1+p^{1-s}}{p+p^{-s}}
     \prod_{p|N/(N,d)}\frac{1+(p-2)p^{-s}}{p-p^{-s}}.
  \end{split}
  \end{equation*}
  By Lemma \ref{lemma: kappa(0)}
  \begin{equation*}
  \begin{split}
    \kappa_0^-(0)&=2\frac{\Lambda'(1,\chi_d)}{\Lambda(1,\chi_d)}
     -\sum_{p|D/(D,d)}\left(\frac{-p^{1-s}\log p}{1+p^{1-s}}
     -\frac{-p^{-s}\log p}{p+p^{-s}}\right)\Big|_{s=0} \\
  &\qquad\qquad-\sum_{p|N/(N,d)}\left(\frac{-(p-2)p^{-s}\log p}{1+(p-2)p^{-s}}
     -\frac{p^{-s}\log p}{p-p^{-s}}\right)\Big|_{s=0} \\
  &=2\frac{\Lambda'(1,\chi_d)}{\Lambda(1,\chi_d)}
   +\sum_{p|D/(D,d)}\frac{p-1}{p+1}\log p
   +\sum_{p|N/(N,d)}\log p
  \end{split}
  \end{equation*}
  and the proof of the lemma is complete.
\end{proof}

\begin{Example} Let $\psi_{F_f}(\tau)$ and $\psi_{F_g}(\tau)$ be the
  Borcherds forms given in Proposition
  \ref{proposition: realization 2}. In this example, we shall utilize
  Lemmas \ref{lemma: kappa(m)} and \ref{lemma: kappa0(0)} to
  determine the absolute values of $\psi_{F_f}(\tau)$ at the CM-point
  of discriminant $-4$ and that of $\psi_{F_g}(\tau)$ at
  the CM-point of discriminant $-3$.

  Let $B=\JS{-1,3}\Q$, $\O=\Z+\Z I+\Z J+\Z IJ$, and the embedding
  $\iota:B\hookrightarrow M(2,\R)$ be chosen as in Section
  \ref{section: Schwarzian}. Choose $\lambda=I$. Then $\phi:i\to I$
  defines an optimal embedding relative to $(\O,\Z[i])$ and the fixed
  point $\tau_d$ of $\iota(\phi(I))$ in the upper half-plane is a
  CM-point of discriminant $d=-4$. By Theorem \ref{theorem: Schofer},
  Lemma \ref{lemma: kappa0(0)}, and \eqref{equation: Lambda'}, we have
  \begin{equation*}
  \begin{split}
    \log\left|\psi_{F_f}(\tau_d)(\Im\tau_d)^4\right|
  &=-\frac14\left(2\kappa_0(3)+8\kappa_0(0)+8\Gamma'(1)+8\log(2\pi)\right)\\
  &=-\frac12\kappa_0(3)
    -4\frac{\Lambda'(1,\chi_d)}{\Lambda(1,\chi_d)}-\log 3
    -2\Gamma'(1)-2\log(2\pi) \\
  &=-\frac12\kappa_0(3)+4\frac{L'(0,\chi_d)}{L(0,\chi_d)}
    -\log 3+2\log|d|-2\log(8\pi^2).
  \end{split}
  \end{equation*}
  The term that needs some work is $\kappa_0(3)$.

  We have $L_+=\Z I$ and $L_-=\Z J+\Z IJ$. Thus, $L=L_++L_-$ and
  by \eqref{equation: kappa}, we have
  $$
    \kappa_0(3)=\sum_{x\in L_+}\kappa_0^-(3-\gen{x,x}/2)
   =\kappa_0^-(3)+2\kappa_0^-(2).
  $$
  With respect to the basis $\{J,IJ\}$, the Gram matrix of $L_-$ is
  $\SM{-6}00{-6}$. Thus, the sets $S_{m,\mu}$ in \eqref{equation: S}
  is $\{2,3\}$ for both $\kappa_0^-(3)$ and $\kappa_0^-(2)$.
  Using results in Section 4 of \cite{Kudla-Yang}, we find that
  $$
    W_{3,2}(s,\varphi_{0,2})=\frac12(1-2^{-2s}), \qquad
    W_{3,3}(s,\varphi_{0,3})=\frac13(1+2\cdot3^{-s}+3^{-2s}),
  $$
  and
  $$
    W_{2,2}(s,\varphi_{0,2})=\frac12(1+2^{-3s}), \qquad
    W_{2,3}(s,\varphi_{0,3})=\frac13(1-3^{-s}).
  $$
  Therefore, by Lemma \ref{lemma: kappa(m)},
  $$
    \kappa_0^-(3)=-8\log 2, \qquad \kappa_0^-(2)=-2\log 3,
  $$
  and $\kappa_0(3)=-8\log 2-4\log 3$. It follows that
  \begin{equation} \label{equation: Ff value}
    \left|\psi_{F_f}(\tau_d)(\Im\tau_d)^4\right|
   =48\cdot\frac{|d|^2}{64\pi^4}e^{4L'(0,\chi_d)/L(0,\chi_d)}.
  \end{equation}

  We next determine the value of $\psi_{F_g}(\tau)$ at the CM-point of
  discriminant $d=-3$. Choose $\lambda=3I-J+IJ$ so that
  $\phi:\sqrt{-3}\to\lambda$ defines an optimal embedding of
  discriminant $-3$. By Theorem \ref{theorem: Schofer}, Lemma
  \ref{lemma: kappa0(0)}, and \eqref{equation: Lambda'} again, we have
  \begin{equation*}
  \begin{split}
    \log\left|\psi_{F_g}(\tau_d)(\Im\tau_d)^6\right|
  &=-\frac12\kappa_0(1)-6\frac{\Lambda'(1,\chi_d)}{\Lambda(1,\chi_d)}
    -\log 2-2\Gamma'(1)-2\log(2\pi) \\
  &=-\frac12\kappa_0(1)+6\frac{L'(0,\chi_d)}{L(0,\chi_d)}
    -\log 2+3\log|d|-3\log(8\pi^2).
  \end{split}
  \end{equation*}

  By Corollary \ref{corollary: lattice}, the lattice $L_-$ has
  discriminant $12$ and its Gram matrix must be equivalent to
  $\SM{-4}{-2}{-2}{-4}$. Since the discriminant of the lattice
  $L_++L_-$ is equal to that of $L$, $L/(L_++L_-)$ is trivial.
  Consequently,
  $$
    \kappa_0(1)=\sum_{x\in L_+}\kappa_0^-(1-\gen{x,x}/2)
   =\kappa_0^-(1).
  $$
  The set $S_{m,\mu}$ in \eqref{equation: S} is $\{2,3\}$ for
  $\kappa_0^-(1)$. Using Theorems 4.3 and 4.4 of \cite{Kudla-Yang}, we
  find
  $$
    W_{1,2}(s,\varphi_{0,2})=\frac12(1-2^{-s}), \qquad
    W_{1,3}(s,\varphi_{0,3})=\frac1{\sqrt3}(1+3^{-s}).
  $$
  Then, Lemma \ref{lemma: kappa(m)} yields
  $$
    \kappa_0^-(1)=-6\sqrt3\cdot\frac{\log 2}{3}\cdot\frac2{\sqrt3}
   =-4\log 2.
  $$
  Finally, we arrive at
  \begin{equation} \label{equation: Fg value}
    \left|\psi_{F_g}(\tau_d)(\Im\tau_d)^6\right|
   =2\cdot\frac{|d|^3}{512\pi^6}e^{6L'(0,\chi_d)/L(0,\chi_d)}.
  \end{equation}
\end{Example}

\begin{Corollary} \label{corollary: C}
  The absolute values of the constants $C_1$ and $C_2$ in
  \eqref{equation: C1} and \eqref{equation: C2} are
  $$
    |C_1|=\frac{12}{\pi^4}e^{4L'(0,\chi_{-4})/L(0,\chi_{-4})}, \qquad
    |C_2|=\frac{27(1+\sqrt3)^6}{256\pi^6}
        e^{6L'(0,\chi_{-3})/L(0,\chi_{-3})},
  $$
  respectively.
\end{Corollary}

\begin{proof} The CM-point of discriminant $-4$ in the
  example above is $\tau_{-4}=i$. According to our choice of $s(\tau)$
  in Proposition \ref{proposition: realization 1}, we have $s(i)=0$.
  Therefore, the right-hand side of \eqref{equation: C1} is simply
  $C_1$. Then \eqref{equation: Ff value} gives us the absolute value
  of $C_1$. The determination of $|C_2|$ is similar.
\end{proof}

\begin{Remark} The values of $|C_1|$ and $|C_2|$ can also be
  determined by considering the values of the Borcherds forms at the
  CM-point $\tau_{-24}$ of discriminant $-24$. At the point
  $\tau_{-24}$, the functions $s(\tau)$ and $t(\tau)$ take value $1$.
  Thus, the right-hand sides of \eqref{equation: C1} and
  \eqref{equation: C2} can be expressed in terms of Gamma values using
  Gauss' formula
  $_{}2F_1(a,b;c;1)=\Gamma(c)\Gamma(c-a-b)/(\Gamma(c-a)\Gamma(c-b))$.
  By repeatedly applying Euler's reflection formula and Gauss'
  multiplication formula, we arrive at the same expressions for
  $|C_1|$ and $|C_2|$.
\end{Remark}



\begin{Example}
  Consider the case $d=-163$. By Theorem \ref{theorem: Schofer} and
  Lemma \ref{lemma: kappa0(0)}, we have
  $$
    \log\left|\psi_{F_f}(\tau_d)(\Im\tau_d)^4\right|
   =-\frac12\kappa_0(3)+4\frac{L'(0,\chi_d)}{L(0,\chi_d)}
    -\log3-\frac23\log 2+2\log|d|-2\log(8\pi^2).
  $$
  On Page 851 of \cite{Errthum}, it is computed that
  $$
    \kappa_0(3)=-\frac{40}3\log 2-4\log 3-4\log 5-4\log(11)-4\log(17).
  $$
  Thus,
  $$
    \left|\psi_{F_f}(\tau_d)(\Im\tau_d)^4\right|
   =2^6\cdot3\cdot5^2\cdot11^2\cdot17^2\cdot\frac{|d|^2}{64\pi^4}
    e^{4L'(0,\chi_d)/L(0,\chi_d)}.
  $$
\end{Example}

We now give the values of the Borcherds forms $\psi_{F_f}(\tau)$ and
$\psi_{F_g}(\tau)$ at various CM-points. The computation is done using
Magma \cite{Magma}. (The use of Magma is not essential. We use Magma only
because it has built-in functions for computation about quaternion
algebras.) The Magma code is available as an accompanying file to this
paper.

\begin{Lemma} \label{lemma: values}
  For a fundamental discriminant $d<0$ appearing in Theorem
  \ref{theorem: main theorem}, let $\tau_d\in\H^+$ be a CM-point of
  discriminant $d$, and
  $$
    \omega_d=e^{L'(0,\chi_d)/2L(0,\chi_d)}
   =\frac1{\sqrt{|d|}}\prod_{a=1}^{|d|-1}\Gamma\left(
    \frac a{|d|}\right)^{\chi_d(a)\mu_d/4h_d}.
  $$
  Let $A_d$ be the number such that
  \begin{equation} \label{equation: Ad}
    \left|\psi_{F_f}(\tau_d)(\Im\tau_d)^4\right|=A_d\frac{|d|^2}{64}\left(
    \frac{\omega_d}{\sqrt\pi}\right)^8.
  \end{equation}
  Then we have
  $$ \extrarowheight3pt
  \begin{array}{||c|l||c|l||c|l||} \hline\hline
  d  & A_d & d & A_d & d & A_d \\ \hline
  -4 & 2^4\cdot 3 & -132 & 2^4\cdot3^2\cdot5^2
     & -148 & 2^4\cdot3\cdot5^2\cdot17^2 \\
 -24 & 2^4\cdot3^2 & -43 & 2^6\cdot3\cdot5^2
     & -232 & 2^4\cdot3\cdot5^2\cdot23^2\cdot29 \\
-120 & 2^4\cdot3^3\cdot5 & -88 & 2^4\cdot3\cdot5^2\cdot 11
     & -708 & 2^4\cdot3^2\cdot5^2\cdot17^2\cdot29^2\\
 -52 & 2^4\cdot3\cdot5^2 & -312 & 2^4\cdot3^2\cdot5^2\cdot11^2
     & -163 & 2^6\cdot3\cdot5^2\cdot11^2\cdot17^2 \\
  \hline\hline
  \end{array}
  $$
  Also, let $B_d$ be the number such that
  \begin{equation} \label{equation: Bd}
    \left|\psi_{F_g}(\tau_d)(\Im\tau_d)^6\right|=B_d
    \frac{|d|^3}{512}\left(\frac{\omega_d}{\sqrt\pi}\right)^{12}.
  \end{equation}
  We have
  $$ \extrarowheight3pt
  \begin{array}{||c|l||c|l||c|l||} \hline\hline
  d & B_d & d & B_d & d & B_d \\ \hline
  -3 & 2 & -19 & 2\cdot3^2 & -67 & 2\cdot3^2\cdot7^2\cdot11^2 \\
 -84 & 2^4\cdot7 & -168 & 2^3\cdot7\cdot11^2 & 
    -372 & 2^4\cdot7^2\cdot19^2\cdot31 \\
 -40 & 2^3\cdot3^2 & -228 & 2^6\cdot7^2\cdot19
     & -408 & 2^3\cdot7^2\cdot11^2\cdot31^2 \\
 -51 & 2\cdot7^2 & -123 & 2\cdot7^2\cdot19^2
     & -267 & 2\cdot7^2\cdot31^2\cdot43^2 \\ \hline\hline
  \end{array}
  $$
\end{Lemma}
\end{section}

\begin{section}{Proof of Theorem \ref{theorem: main theorem}}
  In this section, we shall convert informations from Lemma \ref{lemma:
    values} into special-value formulas for hypergeometric functions.

  We retain our choices of $B$, $\O$, $\iota$, the fundamental
  domain, and etc. from Section \ref{section: Schwarzian}. In the
  following discussion, we let $s$ be the Hauptmodul of $X_0^6(1)/W_6$
  that takes values $0$, $1$, and $\infty$ at the CM-points of
  discriminants $-4$, $-24$, and $-3$, respectively. According to the
  choice of the Fundamental domain in Section \ref{section: Schwarzian}, these
  CM-points are represented by $i$, $(\sqrt6-\sqrt2)i/2$, and
  $(-1+i)/(1+\sqrt3)$, respectively. Let also $t=1/s$.
  For a CM-point $\tau_d$ of a fundamental discriminant $d<0$
  inside the fundamental domain, we let $\phi:\Q(\sqrt
  d)\hookrightarrow B$ be the corresponding optimal embedding and
  assume that $\phi(\sqrt{d})=a_1I+a_2J+a_3IJ$. Then we have
  \begin{equation} \label{equation: tau}
    \tau_d=\frac{a_2\sqrt3+\sqrt{d}}{a_1+a_3\sqrt3}.
  \end{equation}
  We first recall a technical lemma from \cite{Yang-Ramanujan}.

\begin{Lemma}[{\cite[Lemma 5]{Yang-Ramanujan}}] \label{lemma: location}
  If $s(\tau_d)$ takes a value in the line segment $[0,1]$, then
  $a_2=0$. If $s(\tau_d)$ takes a value in $[1,\infty)$, then
  $a_1=3a_3$. If $s(\tau_d)$ takes a negative value, then $a_2=-a_3$.
\end{Lemma}

Recall that $\psi_{F_f}(\tau)$ and $\psi_{F_g}(\tau)$ are the
Borcherds forms defined in \eqref{equation: C1} and \eqref{equation:
  C2}, respectively.

\begin{Proposition}
  \label{proposition: main proposition}
  Assume that $-1<s(\tau_d)<1$. Let $A_d$ be the real
  number such that \eqref{equation: Ad} holds. Then we have
  \begin{equation} \label{equation: F1 value}
    {}_2F_1\left(\frac1{24},\frac5{24};\frac34;s(\tau_d)\right)^8
   =\frac{A_d}{2^{12}\cdot3}(a_1+\sqrt{|d|})^4\left(\frac{\omega_d}{\omega_{-4}}
    \right)^8
  \end{equation}
  and
  \begin{equation} \label{equation: F2 value}
    {}_3F_2\left(\frac13,\frac12,\frac23;\frac34;\frac54;s(\tau_d)\right)^4
   =\frac{3^2A_d}{2^{10}|s(\tau_d)|}(a_2^2+a_3^2)^2\omega_d^8.
  \end{equation}
  Assume that $-1<t(\tau_d)<1$. Let $B_d$ be the real number such that
  \eqref{equation: Bd} holds. Then we have
  \begin{equation} \label{equation: G1 value}
  \begin{split}
    {}_2F_1\left(\frac1{24},\frac7{24};\frac56;t(\tau_d)\right)^{12}
  &=\frac{B_d}{2^7\cdot3^3}\left(\frac{\omega_d}{\omega_{-3}}\right)^{12} \\
  &\qquad  \times \begin{cases}
    ((a_2+2a_3)\sqrt3+\sqrt{|d|})^6, &\text{if }t(\tau_d)>0, \\
    ((a_1-2a_3)\sqrt3+\sqrt{|d|})^6, &\text{if }t(\tau_d)<0, \end{cases}
  \end{split}
  \end{equation}
  and
  \begin{equation} \label{equation: G2 value}
    {}_3F_2\left(\frac14,\frac12,\frac34;\frac56,\frac76;t(\tau_d)\right)^6
   =\frac{B_d}{216|t|}\omega_d^{12}\times
    \begin{cases}
    27(a_2+a_3)^6, &\text{if }t(\tau_d)>0, \\
    (a_1-3a_3)^6, &\text{if }t(\tau_d)<0. \end{cases}
  \end{equation}
\end{Proposition}

\begin{proof} For convenience, set
  $$
   F_1(s)={}_2F_1(1/24,5/24;3/4;s),  \qquad
   F_2(s)={}_2F_1(7/24,11/24;5/4;s),
  $$
  and $s_d=s(\tau_d)$. Note that we have
  $F_1(s)F_2(s)={}_3F_2(1/3,1/2,2/3;3/4,5/4,s)$. Let
  $C=-1/\sqrt[4]{12}\omega_{-4}^2$. By Lemma 5 of
  \cite{Yang-Ramanujan}, we have
  \begin{equation} \label{equation: F2/F1}
    \frac{Cs_d^{1/4}F_2(s_d)}{F_1(s_d)}=\frac{\tau_d-i}{\tau_d+i}.
  \end{equation}
  Combining \eqref{equation: C1}, \eqref{equation: Ad},
  \eqref{equation: tau}, and Corollary \ref{corollary: C}, we find
  \begin{equation*}
  \begin{split}
     A_d\frac{|d|^2}{64}\left(\frac{\omega_d}{\sqrt\pi}\right)^8
   &=\frac{12\omega_{-4}^8|d|^2}{\pi^4(a_1+a_3\sqrt3)^4}
     F_1(s_d)^8\left|1-\frac{\tau_d-i}{\tau_d+i}\right|^8 \\
   &=\frac{12\omega_{-4}^8|d|^2}{\pi^4(a_1+a_3\sqrt3)^4}
     F_1(s_d)^8\left(\frac{2(a_1+a_3\sqrt3)(a_1-\sqrt{|d|})}{3(a_2^2+a_3^2)}
     \right)^4.
  \end{split}
  \end{equation*}
  Simplifying the identity, we get \eqref{equation: F1 value}. To
  prove \eqref{equation: F2 value}, we observe that from
  \eqref{equation: F2/F1} we obtain
  $$
    F_2(s_d)=\sqrt[4]{12}\omega_{-4}^2\frac{F_1(s_d)}{|s_d|^{1/4}}\left|
    \frac{\tau_d-i}{\tau_d+i}\right|
   =\sqrt[4]{12}\omega_{-4}^2\frac{F_1(s_d)}{|s_d|^{1/4}}
    \frac{a_1-\sqrt{|d|}}{\sqrt{3(a_2^2+a_3^2)}}.
  $$
  Combining this with \eqref{equation: F1 value}, we obtain
  \begin{equation*}
  \begin{split}
    F_1(s_d)^8F_2(s_d)^8
  &=2^4\cdot3^2\cdot\omega_{-4}^{16}\frac{F_1(s_d)^{16}}{s_d^2}
    \frac{(a_1-\sqrt{|d|})^8}{3^4(a_2^2+a_3^2)^4} \\
  &=\frac{A_d^2}{2^{20}\cdot3^4\cdot s_d^2}
    \frac{(a_1+\sqrt{|d|})^8(a_1-\sqrt{|d|})^8}{(a_2^2+a_3^2)^4}\omega_d^{16}
   =\frac{3^4A_d^2}{2^{20}s_d^2}(a_2^2+a_3^2)^4\omega_d^{16}.
  \end{split}
  \end{equation*}
  Simplifying the equality, we obtain \eqref{equation: F2 value}.

  Similarly, we write $t_d=t(\tau_d)$, and
  $$
    G_1(t)={}_2F_1\left(\frac1{24},\frac7{24};\frac56;t\right), \qquad
    G_2(t)={}_2F_1\left(\frac5{24},\frac{11}{24};\frac76;t\right).
  $$
  Then $G_1(t)G_2(t)={}_3F_2(1/4,1/2,3/4;5/6,7/6;t)$. Let
  $C'=e^{-2\pi i/8}/\sqrt[6]2\omega_{-3}^2$. We have
  $$
    \frac{C't_d^{1/6}G_2(t_d)}{G_1(t_d)}
   =\frac{\tau_d-\tau_{-3}}{\tau_d-\overline\tau_{-3}}, \qquad
    \tau_{-3}=\frac{-1+i}{1+\sqrt3}.
  $$
  Using
  $$
    \left|\frac{\tau_d-\tau_{-3}}{\tau_d-\overline\tau_{-3}}\right|^2
   =\frac{\sqrt3(a_1+a_2-a_3)-\sqrt{|d|}}{\sqrt3(a_1+a_2-a_3)+\sqrt{|d|}},
  $$
  $$
    \left|1-\frac{\tau_d-\tau_{-3}}{\tau_d-\overline\tau_{-3}}\right|^2
   =\frac2{1+\sqrt3}\frac{a_1+a_3\sqrt3}{\sqrt3(a_1+a_2-a_3)+\sqrt{|d|}},
  $$
  \eqref{equation: C2}, \eqref{equation: Fg value},
  and Corollary \ref{corollary: C}, we deduce that
  \begin{equation*}
  \begin{split}
    B_d\frac{|d|^3}{512}\left(\frac{\omega_d}{\sqrt\pi}\right)^{12}
  &=\frac{27(1+\sqrt3)^6\omega_{-3}^{12}|d|^3}{256\pi^6(a_1+a_3\sqrt3)^6}G_1(t_d)^{12}
    \left|1-\frac{\tau_d-\tau_{-3}}{\tau_d-\overline\tau_{-3}}\right|^{12} \\
  &=\frac{27\omega_{-3}^{12}|d|^3}{4\pi^6(\sqrt3(a_1+a_2-a_3)+\sqrt{|d|})^6}
    G_1(t_d)^{12},
  \end{split}
  \end{equation*}
  so that
  $$
    G_1(\tau_d)^{12}=\frac{B_d(\sqrt3(a_1+a_2-a_3)+\sqrt{|d|})^6}
    {2^7\cdot3^3}\left(\frac{\omega_d}{\omega_{-3}}\right)^{12}
  $$
  and
  \begin{equation*}
  \begin{split}
    G_1(\tau_d)^{12}G_2(\tau_d)^{12}
  &=\frac{4\omega_{-3}^{24}}{t_d^2}G_1(t_d)^{24}\left(
    \frac{\sqrt3(a_1+a_2-a_3)-\sqrt{|d|}}{\sqrt3(a_1+a_2-a_3)+\sqrt{|d|}}
    \right)^6 \\
  &=\frac{B_d^2}{2^{12}\cdot3^6}\left(3(a_1+a_2-a_3)^2-|d|\right)^6\omega_d^{24}\\
  &=\frac{B_d^2}{2^6\cdot3^6}\left(
    a_1^2+3a_2^2+3a_3^2+3a_1a_2-3a_2a_3-3a_1a_3\right)^6\omega_d^{24}.
  \end{split}
  \end{equation*}
  With Lemma \ref{lemma: location}, these two identities reduce to
  \eqref{equation: G1 value} and \eqref{equation: G2 value},
  respectively. This completes the proof.
\end{proof}

\begin{proof}[Proof of Theorem \ref{theorem: main theorem}]
The values of $s(\tau)$ and $t(\tau)$ at CM-points were computed in
\cite{Errthum}. They are the rational numbers $M/N$ from the two
tables in Theorem \ref{theorem: main theorem}. The optimal embeddings
corresponding to the CM-points inside the fundamental domain are given
in the two tables below.
$$ \extrarowheight3pt
\begin{array}{||r|l||r|l||} \hline\hline
 d & \phi(\sqrt d) & d & \phi(\sqrt d) \\ \hline
-52 & 8I+2IJ & -120 & 12I-2J+2IJ \\
-88 & 10I+2IJ & -43 & 7I-J+IJ\\
-132& 12I+2IJ & -232 & 16I-2J+2IJ\\
-312& 18I+2IJ & -163 & 13I-J+IJ\\
-148& 14I+4IJ & &\\
-708& 30I+8IJ & &\\
\hline\hline
\end{array}
$$
$$ \extrarowheight3pt
\begin{array}{||r|l||r|l||} \hline\hline
 d & \phi(\sqrt d) & d & \phi(\sqrt d) \\ \hline
-84 & 12I-2J+4IJ & -40 & 8I-2J+2IJ \\
-51 & 9I-J+3IJ & -19 & 5I-J+IJ \\
-168&18I-4J+6IJ & -228 & 18I-4J+4IJ \\
-123&15I-3J+5IJ & -67 & 11I-3J+3IJ \\
-372&24I-2J+8IJ & &\\
-408&30I-8J+10IJ& &\\
-267&21I-3J+7IJ & & \\
\hline\hline
\end{array}
$$
Here the left columns of the two tables are for discriminants $d$ with
$s(\tau_d)>0$ and $t(\tau_d)>0$, respectively.
Combining informations from Lemma \ref{lemma: values}, Proposition
\ref{proposition: main proposition}, and the above two tables, we
obtain the identities in Theorem \ref{theorem: main theorem}.
\end{proof}
\end{section}

\begin{section}{Further examples}
\label{section: additional}
Observe that for each discriminant $d$ appearing in Theorem
\ref{theorem: main theorem}, there is only one CM-point of
discriminant $d$ on the Shimura curve $X_0^6(1)/W_6$. In such cases,
Schofer's formula readily tells us the absolute value of a
Borcherds form at the unique CM-point of discriminant $d$. However,
in general, we can only read from Schofer's formula the products of
values of Borcherds forms at CM-points. In this section, we
introduce a technique to separate the value at a CM-point from those
at the other CM-points of the same discriminant using Hecke operators.
This technique relies on the method developed in
\cite{Yang-Schwarzian} for computing Hecke operators. Here we will
work out the case $d=-276$. In principle, the method works at least
for any imaginary quadratic number field whose ideal class
group, after quotient by the prime ideals lying above $2$ and $3$, is
an elementary $2$-group.

Let $E=\Q(\sqrt{-276})$ and $R$ be the ring of integers in $E$. There
are two CM-points of discriminant $d=-276$ on $X_0^6(1)/W_6$,
represented by the two points
$$
  \tau_1=\frac{\sqrt{-69}}{9+2\sqrt3}, \qquad
  \tau_2=\frac{-3\sqrt3+\sqrt{-69}}{12+4\sqrt3}
$$
in the fundamental domain. The corresponding optimal embeddings
$\phi_1$ and $\phi_2$ are
$$
  \lambda_1=\phi_1(\sqrt{-276})=18I+4IJ, \qquad
  \lambda_2=\phi_2(\sqrt{-276})=24I-6J+8IJ,
$$
respectively. According to the table at the end of Section 5 of
\cite{Yang-ModEqs}, the values of the Hauptmodul $s(\tau)$ at these two
points are $(166139596\pm 95538528\sqrt3)/1771561$. (The values can
also be determined using Borcherds forms and Schofer's formula.) From
Lemma \ref{lemma: location}, we deduce that
$$
  s(\tau_1)=\frac{166139596-95538528\sqrt3}{1771561}, \qquad
  s(\tau_2)=\frac{166139596+95538528\sqrt3}{1771561}.
$$
Call these two numbers $s_1$ and $s_2$, respectively. Let $\fp_2$ and
$\fp_3$ be the prime ideals of $R$ lying above $2$ and $3$,
respectively, and let $\fp_5$ be any prime above $5$. Then the ideal
class group of $R$ is isomorphic to $(\Z/2\/Z)\times(\Z/4\Z)$
generated by the element $\fp_2$ of order $2$ and the element $\fp_5$
of order $4$. Moreover, the product $\fp_2\fp_3\fp_5^2$ is a principal
ideal. It follows that the ideal class group, afte quotient by the
subgroup generated by $\fp_2$ and $\fp_3$, is cyclic of order $2$ and
generated by $\fp_5$. In terms of CM-points on $X_0^6(1)/W_6$, this
means that there should exist an element $\alpha$ of norm $5$, $10$, $15$,
or $30$ in $\O$ such that $\iota(\alpha)\tau_1=\tau_2$. (Here we
retain the notations $\O$, $\iota$, and etc. used in Section
\ref{section: Schwarzian}.) Indeed,
such an element is
$$
  \alpha=3-2I-IJ.
$$
(Another element is $\alpha'=(3-9I+J-3IJ)/2$.)
In other words, we have $\lambda_2=\alpha\lambda_1\alpha^{-1}$.

Now let $F(\tau)=\psi_{F_f}(\tau)$ be the modular form of weight $8$ 
defined in Proposition \ref{proposition: realization 2} and set
$$
  \wt F(\tau):=F\big|_8\iota(\alpha)
 =\frac{10^4}{((2+\sqrt3)\tau-3)^8}F\left(\frac{3\tau+2-\sqrt3}
  {(-2-\sqrt3)\tau+3}\right).
$$
In general, we have
$$
  \frac{10^4}{|(2+\sqrt3)\tau-3|^8}
 =\left(\frac{\Im\iota(\alpha)\tau}{\Im\tau}\right)^4.
$$
Thus,
\begin{equation} \label{equation: additional 1}
  \left|F(\tau_2)\right|
 =\left(\frac{\Im\tau_1}{\Im\tau_2}\right)^4\left|\wt F(\tau_1)\right|.
\end{equation}
On the other hand, Schofer's formula yields
$$
  \left|F(\tau_1)F(\tau_2)\right|
  (\Im\tau_1)^4(\Im\tau_2)^4
 =2^8\cdot3^4\cdot11^2\left(\frac{|d|^2}{64\pi^4}\omega_d^8\right)^2.
$$
Substituting \eqref{equation: additional 1} into this, we obtain
\begin{equation} \label{equation: additional 2}
  \left|F(\tau_1)(\Im\tau_1)^4\right|^2\left|
  \frac{\wt F(\tau_1)}{F(\tau_1)}\right|
 =2^8\cdot3^4\cdot11^2\left(\frac{|d|^2}{64\pi^4}\omega_d^8\right)^2.
\end{equation}
The main task remained is to determine the value of $\wt
F(\tau_1)/F(\tau_1)$.

Let $\Gamma$ be the discrete subgroup of $\mathrm{PSL}(2,\R)$ such
that $X_0^6(1)/W_6=\Gamma\backslash\H^+$, i.e.,
$\Gamma:=\{\iota(\gamma)/(\det\gamma)^{1/2}:\gamma\in N_B^+(\O)\}$.
Let $\gamma_j$, $j=1,\ldots,5$, be elements in
$\Gamma\iota(\alpha)\Gamma$ such that $\gamma_0=\iota(\alpha)$ and
$\gamma_j$, $j=1,\ldots,5$, form a complete set of coset
representatives of $\Gamma\backslash\Gamma\iota(\alpha)\Gamma$.
In Section 4 of \cite{Yang-Ramanujan}, by using results from
\cite{Yang-Schwarzian}, we find that
\begin{equation*}
\begin{split}
  \prod_{j=0}^5\left(y-\frac{F\big|_8\gamma_j}{F}\right)
&=y^6+\frac{114}{125}y^5-\frac{6333}{78125}y^4
 +\frac{4}{5^{11}}(8640000s-5177953)y^3 \\
&\qquad+\frac{3}{5^{15}}(8467200000s+1804020097)y^2 \\
&\qquad+\frac{726}{5^{20}}(93744000000s-3501556201)y \\
&\qquad+\frac1{5^{16}}(138240s+14641)^2.
\end{split}
\end{equation*}
Substituting $s$ by $s_1=(166139596-95538528\sqrt3)/1771561$, we
deduce that $\wt F(\tau_1)/F(\tau_1)$ is a zero of
\begin{equation*}
\begin{split}
  (9150625y^2+(40464094y-20903960\sqrt3)y+82650625-47425000\sqrt3)g(y),
\end{split}
\end{equation*}
where $g(y)\in\Q(\sqrt3)[y]$ is an irreducible polynomial of degree
$4$ over $\Q(\sqrt3)$. In fact, we can show that it is a zero of the
factor of degree $2$ shown above. Hence, we have
\begin{equation} \label{equation: additional 3}
  \left|\frac{\wt F(\tau_1)}{F(\tau_1)}\right|
 =\left(\frac{82650625-47425000\sqrt3}{9150625}\right)^{1/2}
 =\left(\frac{14-5\sqrt3}{11}\right)^2.
\end{equation}
(It is possible to determine the precise value, not just the
absolute value. The two zeros of the
factor of degree $2$ are $\wt F(\tau_1)/F(\tau_1)$ and the value of
$(F\big|_8\iota(\alpha'))/F$ at $\tau_1$, where $\alpha'=(3-9I+J-3IJ)/2$.
It is easy to find the
ratio of the two values and hence determine $\wt F(\tau_1)/F(\tau_1)$.)
Substituting \eqref{equation: additional 3}
into \eqref{equation: additional 2}, we obtain
$$
  \left|F(\tau_1)(\Im\tau_1)^4\right|=144(14+5\sqrt3)
  \frac{|d|^2}{64\pi^4}\omega_d^8.
$$
By Proposition \ref{proposition: main proposition}, this implies that
$$
  {}_2F_1\left(\frac1{24},\frac5{24};\frac34;s_1\right)^8
 =\frac{3(14+5\sqrt3)}{16}(9+\sqrt{69})^4\left(\frac{\omega_{-276}}{\omega_{-4}}
  \right)^8
$$
and
$$
  {}_3F_2\left(\frac13,\frac12,\frac23;\frac34,\frac54;s_1\right)^4
 =\left(\frac{3(16+23\sqrt3)}{11}\right)^4\left(\frac{2+3\sqrt3}{23}\right)^2
  (2+\sqrt3)\omega_{-276}^8.
$$
\end{section}

\end{document}